\documentclass[12pt]{iopart}

\usepackage{subfig}
\usepackage{iopams} 
\usepackage{graphicx}
\expandafter\let\csname equation*\endcsname\relax
\expandafter\let\csname endequation*\endcsname\relax

\usepackage{lipsum}
\usepackage{graphicx}
\usepackage{epstopdf}
\ifpdf
  \DeclareGraphicsExtensions{.eps,.pdf,.png,.jpg}
\else
  \DeclareGraphicsExtensions{.eps}
\fi

\usepackage{float}
\usepackage{bm}
\usepackage{multirow} 
\usepackage{tabularx}
\usepackage{gensymb}
\usepackage{algpseudocode}

\usepackage{mathtools}
\usepackage[section]{placeins}
\newcommand{\R}{\mathbb{R}}

\DeclarePairedDelimiterX{\norm}[1]{\lVert}{\rVert}{#1}
\newcommand{\palentir}{PaLEnTIR}
\newcommand{\Rbb}{\mathbb{R}}

\newcommand{\Amb}{\mathbf{A}}   
 
\newcommand{\Cmb}{\mathbf{C}} 
 \newcommand{\dmb}{\mathbf{d}}
 
 \newcommand{\fmb}{\mathbf{f}}

 \newcommand{\pemb}{\mathbf{p}}
 
\newcommand{\Rmb}{\mathbf{R}} \newcommand{\rmb}{\mathbf{r}}

 \newcommand{\wmb}{\mathbf{w}}

\newcommand{\betamb}{\boldsymbol{\beta}}
\newcommand{\gammamb}{\boldsymbol{\gamma}}

\newcommand{\chimb}{\boldsymbol{\chi}}

\usepackage{upgreek}
\newcommand{\bfp}{\mathbf{p}}
\newcommand{\bfd}{\mathbf{d}}
\newcommand{\bfr}{\mathbf{r}}
\newcommand{\bfpp}{ {\mathbf{f}}({\mathbf{p}}) }

\begin{document}

\title{Parametric Level-sets Enhanced To Improve Reconstruction (\palentir{})}

\author{Ege Ozsar$^1$, Misha Kilmer$^2$, Eric de Sturler$^3$, Arvind K. Saibaba$^4$, Eric Miller$^1$
\footnote{This research was supported by the U.S. National Science Foundation under awards 1720291, 1935555,  1934553, 1720398, and 1720305.}}
\begin{indented}
\item[]    $^1$Department of Electrical and Computer Engineering, Tufts University, Medford, MA, USA\\%
    $^2$Department of Mathematics, Tufts University, Medford, MA, USA\\%
    $^3$Department of Mathematics, Virginia Tech, Blacksburg, VA, USA\\%
    $^4$ Department of Mathematics, North Carolina State University, Raleigh, NC, USA\\[2ex]%
\end{indented}


\begin{abstract}
\label{sec:abstract}
We introduce \palentir{}, a significantly enhanced parametric level-set (PaLS) method addressing the restoration and reconstruction of piecewise constant objects. Our key contribution involves a unique PaLS formulation utilizing a single level-set function to restore scenes containing multi-contrast piecewise-constant objects without requiring knowledge of the number of objects or their contrasts. Unlike standard PaLS methods employing radial basis functions (RBFs), our model integrates anisotropic basis functions (ABFs), thereby expanding its capacity to represent a wider class of shapes. Furthermore, \palentir{} improves the conditioning of the Jacobian matrix, required as part of the parameter identification process, and consequently accelerates optimization methods. We validate \palentir{}'s efficacy through diverse experiments encompassing sparse and limited angle of view X-ray computed tomography (2D and 3D), nonlinear diffuse optical tomography (DOT), denoising, and deconvolution tasks using both real and simulated data sets.
\end{abstract}

%
\vspace{2pc}
\noindent{\it Keywords}:   level-set, PaLS, parametric, reconstruction, tomography, piecewise constant

%
%
%

\section{Introduction}
\label{sec:intro}
Inverse problems are pivotal in a broad range of science and engineering applications. In the pursuit of extracting the unknown composition and structure of a medium from indirect observations governed by physical models, researchers often focus on characterizing ``regions of interest" (ROIs). These ROIs may include features such as cancerous tumors in diffuse optical data \cite{arridge}, subsurface contaminants in hydrological data \cite{fakhreddine2016imaging}, or buried objects in electromagnetic data \cite{el2005imaging}. Conventionally, these problems typically involve a computationally intensive image formation step followed by ROI identification \cite{rekanos99}. For problems where data are limited, the initial image formation stage will require potentially complex regularization methods to overcome ill-posedness. Alternatively, direct estimation of ROI geometry and contrasts from data offers a more efficient approach \cite{Kirsch98}. These \textit{shape-based} methods are usually better-posed compared with pixel-based problems; however, topologically complicated shapes can lead to challenges \cite{santosa96}. For example, approaches based on using parametric shapes (circles, ellipsoids, etc.) to describe object geometry require that the number of components for the shape is known 
\textit{a priori} or somehow estimated. For this and many other reasons \cite{aghasi10}, level-set methods have found great use for shape-based inverse problems because of their ability to naturally recover objects whose topology (number of connected components) is not known \textit{a priori}. 

Initially proposed by Osher and Sethian \cite{osher88} for curve propagation modeling, level-set methods subsequently gained traction in solving inverse problems, as pioneered by Santosa \cite{santosa96}. Notable efforts have followed, employing level-set evolution techniques \cite{dorn06,feng03}, particularly in image processing \cite{chan01}. For ill-posed problems, regularization of the level-set function becomes essential. Various strategies have emerged to address this, encompassing pixel-based approaches \cite{doel06}, geometric constraints \cite{osher01}, and finite-dimensional basis function spaces like the parametric level-set function (PaLS) concept. This latter approach, introduced for inverse problems by \cite{aghasi10}, is the basis for the work in this paper.  

The PaLS model has proven to be capable of capturing the topological advantages of a level-set function while avoiding difficulties such as the need for explicit regularization and reinitialization that occur frequently when using traditional level-set methods for inverse problems. Moreover, it was shown in \cite{aghasi10} that the low order representation of the inverse problem makes it possible to use Newton and quasi-Newton methods for determining the PaLS parameters. In recent works, the PaLS model in \cite{aghasi10} and variants have been used across a range of application areas and imaging modalities including geophysics \cite{mcm16}, seismology \cite{7784771} and reservoir monitoring \cite{hov17}, image segmentation \cite{mesadi16}, tomography problems \cite{naik14,NIU2022117819,10111292},  electromagnetic imaging \cite{id20},  and multi-modal imaging \cite{est19}.  

Despite the advantages of PaLS, room for improvement remains. Notably, existing PaLS models and most other level-set methods struggle to recover distributions of objects with multiple contrast values. Approaches like colour \cite{hd21}, vector \cite{zhao96}, and binary level-sets \cite{lie06} rely on a number of level-sets proportionate to the number of objects \cite{mesadi16}, leading to scalability challenges as the object count increases. Furthermore, the choice of basis functions for PaLS raises limitations. Traditionally, a weighted superposition of predetermined basis functions, predominantly radial basis functions (RBFs) \cite{lsd20, semerci12, 6192289}, has been employed. However, RBFs possess only circular cross sections and thus can limit the range of objects that can be represented by the model and a given number of RBFs. Moreover, despite their efficacy, many prevailing PaLS models grapple with numerical conditioning issues. Models employing RBFs can encounter non-uniqueness in parameter representation due to infinitely many parameter pairs yielding the same circular cross sections. As we illustrate empirically in this paper, this  can yield ill-conditioned, even singular, Jacobian matrices — detrimental for quasi-Newton and even for trust-region methods, the performance of which degrades under such conditions \cite{aghasi10}.

Here, we develop Parametric Level-Sets Enhanced to Improve Reconstruction (\palentir) for piecewise constant image reconstruction to address the issues identified in the above paragraph. The specific contributions are as follows.  First, we replace the binary Heaviside function used in the existing PaLS models with a smooth transition function. The resulting model can capture \textbf{multiple unknown contrasts with only a single level-set function}. Thus, the number of parameters to be estimated is independent of the number of contrasts.  To the best of our knowledge, \palentir{} is the first level-set model with this property. The efficacy of this multi-contrast feature is contingent upon its capacity to specify space-varying bounds on the
contrasts at relatively coarse scales. The new model achieves this by parametrically characterizing the spatially varying contrast limits on a sparse grid of points. Second, we replace the RBFs with anisotropic basis functions that produce rotated ellipsoidal cross-sections.  This choice provides {greater geometric flexibility} and enhanced \textbf{shape expressiveness}. The model captures more detail while employing fewer basis functions, proving particularly advantageous for challenging scenarios like the representation of long, flat objects where RBFs encounter difficulties. Finally, we demonstrate empirically that  bounding the expansion coefficients associated with the ABFs  and fixing the centers of the basis functions significantly \textbf{improves the numerical performance} of the method, as evidenced by the reduction of the condition number of the Jacobian matrices generated over the course of the reconstruction process.

The organization of the paper is as follows. In section \ref{sec:formulation}, we define our problem of interest and review the parametric level-set method. 
We also briefly discuss the nonlinear least-squares minimization algorithm,
\texttt{TREGS}, used in our experiments,
at the end of this section.
In section \ref{sec:palentir}, \palentir{} formulations for both 2D and 3D inverse problems are detailed. We also discuss the benefits of the new formulation relative to existing PaLS models. We demonstrate the benefits through a sparse view  2D X-ray computed tomography example. We provide experimental results for \palentir{} in section \ref{sec:experiment}. Specifically, in subsection \ref{sec:experiment_6}, we show the results for a 2D linear problem, namely deconvolution. In subsection \ref{sec:experiment_7}, we demonstrate the performance of \palentir{} on limited angle view and sparse angle multi-contrast 2D inversion of the Radon transforms as well as a 3D limited angle view parallel beam tomography experiment. To test our approach on a nonlinear problem, in subsection \ref{sec:experiment_10}, we show the experimental results of the new approach on diffuse optical tomography (DOT), a severely ill-posed inverse problem. Across this range of experiments, we use both real and synthetic data to show the robustness of our approach in a range of data-limited applications. Conclusions are provided in section \ref{sec:conclusion}.

\section{Problem formulation}
\label{sec:formulation}

\subsection{Forward and Inverse Problems}
\label{sec:formulation_1}
Consider a region of space to be imaged, $\Omega \subset \R^d$.
For $\mathbf{r}$ a point in $\Omega$, let us define a spatially-dependent property $f(\mathbf{r})$ of the medium (e.g., electrical conductivity, optical absorption, sound speed, etc.).  We denote with $\mathcal{M}$ the map which takes $f(\mathbf{r})$ to a vector of noise free data.  Typically, $\mathcal{M}$ encompasses the physics of the sensing modality and the engineering details of the associated sensors.  
The data available for processing is equal to $\mathcal{M}(f)$ corrupted by noise.  Unless  specified otherwise, we assume that the noise is additive Gaussian.  
In summary, the forward model can be written 
\begin{equation}
    \mathbf{d}= \mathcal{M}(f) + {\bf w},
    \label{forward}
\end{equation}
where ${\bf w,d} \in \mathbb{R}^{N_{pts}}$ represent the  additive noise and the data available for processing.  

The inverse problem requires determination of the unknown function $f$ from $\mathbf{d}$. Following the penalized likelihood approach, \eqref{forward} leads to the following minimization problem as the basis for recovering $f(\rmb)$,
\begin{equation}
    \label{eq:ip1}
    \min_{f} \frac{1}{2} \|  \mathcal{M}(f)-\dmb \|_2^2 + \xi (f).
\end{equation}
The first term in \eqref{eq:ip1} arises from a Gaussian assumption on the noise and quantifies the mismatch error between model prediction and the data, while $\xi(f)$ is the regularization functional which is usually used when the problem is ill-posed and is chosen based on prior knowledge concerning, e.g., the degree of smoothness associated with $f$ \cite{tikhonov77}. It is also possible to regularize the problem using a geometric parameterization of the unknown property. This is done by either embedding the regularization implicitly in the parameterization of the unknown property, in which case no explicit $\xi$ may be needed, or by expressing it as geometric constraints on the shape of the unknown \cite{dorn06}. In this paper, we follow the former approach using a PaLS type of model.

\subsection{Parametric Level-set methods} 
\label{sec:formulation_2}
Under a PaLS model, 
$f : \mathbb{R}^d \rightarrow \mathbb{R}$
consists of object, $O$, and background, $\Omega \backslash O$, terms and is written as 
\begin{equation}
    \label{eq:f1}
    f(\rmb)= f_O(\rmb)\chi_O(\rmb) + f_B(\rmb)(1-\chi_O(\rmb)) ,
\end{equation}
where $f_O(\rmb)$ and $f_B(\rmb)$ are the generally spatially dependent property values of the object and background, respectively, and $\chi_O(\rmb) = 1 \text{ for }\rmb \in O \text{ and } 0 \text{ elsewhere}$ is the characteristic function of the region $O$.
Under this model, we seek $O$ (or the boundary, $\partial O$), $f_O(\rmb)$, and $f_B(\rmb)$.  
As $O$ can be multiply connected with each component having no specific (i.e., easily parameterized) shape, level-set methods have proven convenient \cite{aghasi10}. Mathematically, the level-set representation of  $O$ satisfies
\begin{equation}
    \phi_O(\rmb) 
        \begin{cases}
            > c & \rmb \in O , \\
            = c & \rmb \in \partial O , \\
            < c & \rmb \in \Omega \backslash O ,
        \end{cases}
\end{equation}
where $c$ is a constant determining the level-set\footnote{This could be taken to be the zero level-set, though in \cite{aghasi10}, the case was made for $c$ slightly larger.} \cite{aghasi10}. In terms of $\phi_O$,  $\chi_O(\rmb) = H(\phi_O(\rmb)-c)$ where $H(x)$ is the Heaviside function.  Thus, \eqref{eq:f1} becomes 
\begin{equation}
    \label{property_level}
    f(\rmb) = f_OH(\phi_O(\rmb)-c)+f_B(1-H(\phi_O(\rmb)-c))
\end{equation}
where $f_O$ and $f_B$ are, for now, taken to be constants.  

Many level-set methods follow a finite difference discretization of the level-set function, which requires a dense collection of nodes. The difficulty of implementing this approach as well as the numerical considerations of its discrete computation overshadow the advantages of the level-set function, especially in the case of ill-posed inverse problems \cite{aghasi10,ben07}. Alternatively, a parametric form for the level-set function expands $\phi_O(\rmb)$ as a linear superposition of a set of basis functions (e.g., polynomial, radial basis function, trigonometric, etc).  Specifically, the original PaLS model takes the form
\begin{equation}
        \phi_{\rm rbf}(\rmb; {\bf p} )= \sum^N_{j=1}\alpha_j \psi(\beta_j(\rmb-\boldsymbol{\chi}_{j})), \  {\bf p}:= \begin{bmatrix} 
            \boldsymbol{\alpha} \\ \boldsymbol{\beta} \\ \boldsymbol{\chi}
        \end{bmatrix} .
    \label{phi_og}
\end{equation}
In this case, the PaLS function is formed as a weighted summation of $N$ basis functions $\psi_j(\rmb):= \psi(\beta_j(\rmb-\boldsymbol{\chi}_{j}))$ for $j=1,2,\dots N.$ The basis functions $\psi_j(\rmb)$ are often taken to be RBFs \cite{9112680,li20,8049349,NIU2022117819}. We will refer to use of such basis functions as ``RBF PaLS'', and it is against such representations that we compare our new \palentir{} representation.\footnote{
In \cite{aghasi10}, the norm in the argument is replaced by a pseudo-norm and $\psi( \cdot )$ is taken to be a compactly supported radial basis function. The compactness of the basis functions may be advantageous in terms of yielding a sparse Jacobian \cite{aghasi10}.  Here, we will use Gaussian RBFs instead of CSRBFs and forgo compactness in favor of retaining a norm, rather than a pseudo-norm.  This choice facilitates the comparative analysis of the expressiveness and numerical properties of the proposed approach over RBF PaLS. The parameterization developed here and the analysis, can be utilized in a CSRBF framework.} Each basis function $\psi_j(\rmb)$ is associated with its own dilation coefficient $\beta_j$, and center location $\boldsymbol{\chi}_{j}$. Ignoring the contrasts $f_B$ and $f_O$ for a moment, $\pemb$ contains the geometric unknowns to be determined, including  $\boldsymbol{\alpha}, \boldsymbol{\beta}$, and $\boldsymbol{\chi}$, containing the expansion coefficients $\alpha_j$, the dilation values $\beta_j$, and PaLS centers $\boldsymbol{\chi}_j\in \Rbb^d$, respectively. In general, the length of this vector is $(d+2)N$, which results in length $4N$ (for 2D problems) or $5N$ (for 3D).

Combining \eqref{phi_og} with \eqref{property_level}, we obtain $f = f(\bfr;\bfp)$.  The goal of a PaLS-based inverse problem is to recover the unknown $\pemb$ based on the observed data $\dmb$ and the model $\mathcal M$. In practice, this requires replacing the exact Heaviside function with a differentiable approximation (we discuss this more in section \ref{sec:palentir_2}). If the contrast coefficients $f_O$ and $f_B$ are known, the resultant inverse problem formulation using PaLS recovers $f$ as $f(\hat{\pemb})$ where $\pemb$ is a minimizer of a nonlinear least squares problem with the objective function
\begin{equation}
    \text{F}(\pemb)= \frac{1}{2} \|  \mathcal{M}(f( {\bf p}))-\dmb \|_2^2.
    \label{optimization}
\end{equation}
When the contrasts are not known, a cyclic descent method is often employed \cite{dorn07}, in which one alternately estimates $\pemb$ for the current estimates of $f_O$ and $f_B$ and then updates the contrasts using the just-computed PaLS parameters.  For linear inverse problems, computing the contrasts is a linear least squares problem which can be solved efficiently, and the bulk of the computational effort is focused on estimating $\pemb$. 

Since the number of PaLS basis functions is always much smaller than the size of the grid obtained from discretizing $\Omega$, the PaLS parameterization acts as a form of regularization, so we can drop the regularization term $\xi(f)$ from the objective function.


In this paper, we use the \texttt{TREGS} algorithm to minimize \eqref{optimization} and solve for the PaLS parameters \cite{TREGS}. The \texttt{TREGS} algorithm is designed to address nonlinear least squares problems encountered in parameterized imaging scenarios. By analyzing spectral components of the Gauss–Newton direction, \texttt{TREGS} determines which components to discard or dampen, minimizing the total number of function and Jacobian evaluations. Leveraging the Basic Trust Region Algorithm \cite{conn_trust_region}, \texttt{TREGS} guarantees global convergence to a critical point, making it a robust optimization tool.   The algorithm employs a discrepancy principle–based stopping criterion. The discrepancy principle states that one should stop iterating when the norm of the residual reaches the norm of the (weighted) noise vector \cite{Hansen_TREGS}. In our experimental setup, an additional stopping criterion is implemented to prevent unnecessary computational expenditure. Specifically, the algorithm concludes its iterations when the relative decrease in the residual falls below a predefined threshold, ensuring judicious use of computational resources.

\section{Parametric Level-sets Enhanced To Improve Reconstruction}  
\label{sec:palentir}
In this section, we introduce \palentir{}, and we discuss its advantages over the RBF PaLS representation. The proposed  \palentir{} model is defined as follows
\begin{align}
    f(\rmb;\pemb)  =&
    C_{H}(\rmb)T_{}\left(\phi(\rmb;\pemb) \right)+  C_{L}(\rmb) \left(1 - T_{}\left(\phi \left(\rmb;\pemb \right) \right)\right) ,
                    \label{eq:newpals}
\\
    \phi(\rmb ; \pemb)  =&\sum^N_{j=1}  \sigma_h \left(\alpha_i \right) 
            \psi \left( \Rmb_j (\rmb-\boldsymbol{\chi}_{j})\right), \     {\bf p}:= \begin{bmatrix} \boldsymbol{\alpha} \\  \boldsymbol{\beta} \\ \boldsymbol{\gamma} \end{bmatrix}.
    \label{eq:newphi}
\end{align}

As explained below, the $N$ matrices, $\Rmb_j := \Rmb({\betamb_j},\gammamb_j)$, for $j=1,2 \dots N$, depend on the $\betamb$ and $\gammamb$ sub-vectors of $\pemb$, respectively. The length of vector $\pemb$ is $3N$ (for 2D) or $7N$ (for 3D).  $T$ is the new transition function, replacing the Heaviside
function,  $C_{H}(\rmb)$ and  $C_{H}(\rmb)$ are the new ``contrast coefficients'', replacing $f_O$ and $f_B$ in \eqref{property_level}, and $\sigma(.)$ is a function to bound the values of $\alpha_i$. Compared to \eqref{phi_og}, the \palentir{} model in \eqref{eq:newpals} enhances PaLS in the following ways, each of which is discussed in depth in the following subsections:
\begin{enumerate}
    \item[A)] \textbf{Multi-contrast, single level-set} reconstructions are obtained by: (a) Replacing the Heaviside function with the smooth transition function $T_{w} :\mathbb{R} \rightarrow (0,1)$; and (b) Relabeling and reinterpreting $f_O(\rmb)$ and $f_B(\rmb)$ as upper and lower contrast bounds, $C_{H}$ and $C_{L}$.  In the simplest case, where there is a single object in the field, these are constant.  For the problems of interest in this paper, both $C_L$ and $C_H$ are functions of space, whose structures and parameterization are discussed in subsection \ref{sec:palentir_2}.
        
    \item[B)] \textbf{Shape-expressiveness} is expanded by replacing the scalar dilation coefficient, $\beta_j$ in PaLS, by a matrix $\Rmb_j$, implementing what we call stretching and sliding.  For 2D problems, each $\Rmb_j$ depends on one element of $\betamb$ and one of $\gammamb$ for a total of $2$ parameters per basis function. In the 3D case, this changes to three $\beta$'s and three $\gamma$'s per $\Rmb_j$, for a total of $6$ parameters per basis function. 
    
    \item[C)] \textbf{Numerical performance} is improved by (a) Fixing the  basis function centers on a grid of pre-specified points, $\chimb_j$, so that these quantities are no longer estimated as part of the inversion process, and combined with the new parameterization, reducing the dimension of the search space;
    (b) Constraining the size of the expansion coefficients using a function $\sigma_h(\alpha_i) \in (-1,1)$; and (c) Replacing RBFs, that have non-unique mapping of parameters to $c$-level-sets, with the new anisotropic basis functions (ABFs).
    \end{enumerate}
We next discuss  each of these three enhancements in turn.  

\subsection{Multi-contrast, Single Level-set Reconstructions}
 \label{sec:palentir_2}

\begin{figure}[!ht]
     \centering
     \subfloat[]
        {\includegraphics[width=0.49\textwidth]{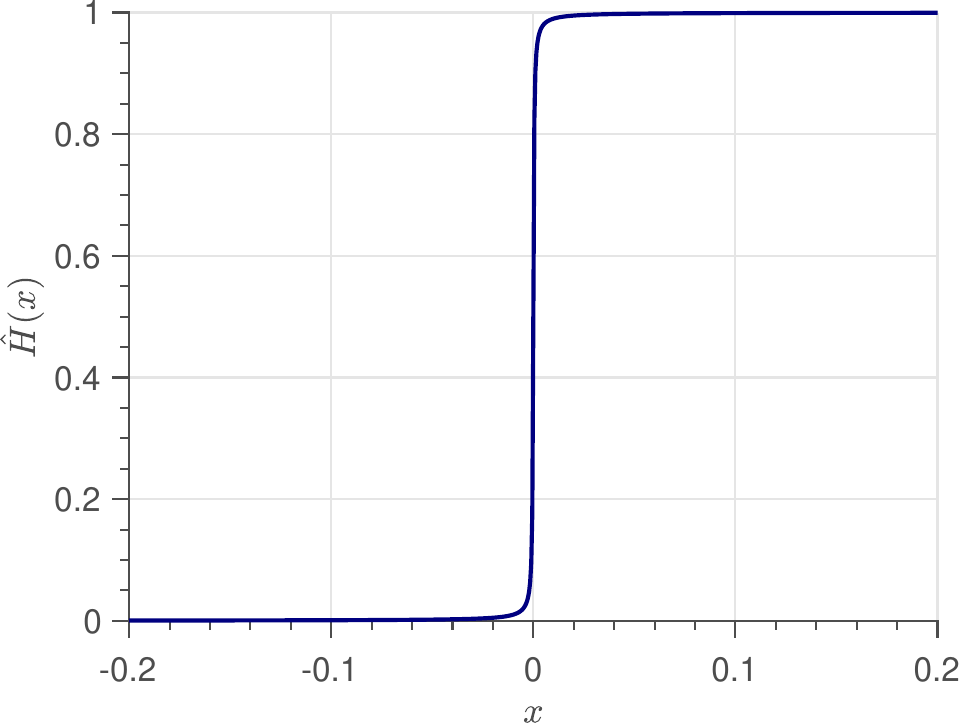}        \label{fig:transition_a}}
        \centering
    \subfloat[]
        {\includegraphics[width=0.49\textwidth]{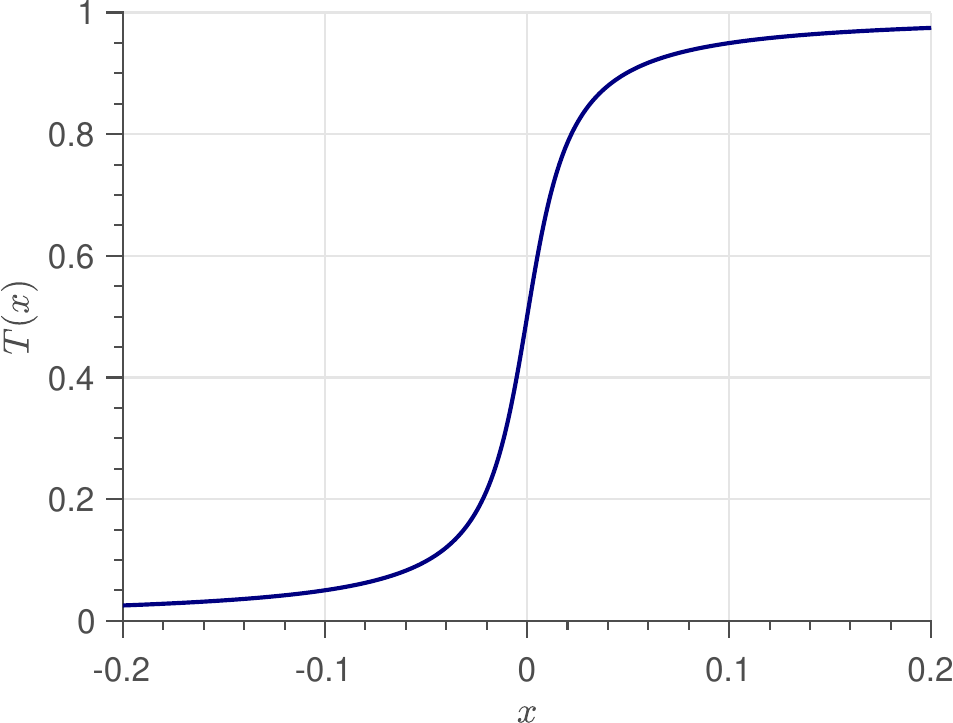}
        \label{fig:transition_b}}

        \centering
    \caption{(a) The plot of the approximate Heaviside function $\hat{H}(x)$ for the zero level-set ($c=0$). (b) The plot of the new transition function $T(x)$ for zero level-set.}
    \label{fig:transition}
\end{figure}

In the previous section, we assumed that the property function $f(\rmb)$ had a binary structure, i.e., each point in $\Omega$ either belongs to the region $O$ or the background. However, it may well be the case that we have to represent regions containing objects with more than two contrast values \cite{shi20, wang04, allaire14}.  To solve this problem, we replace the Heaviside function in the original PaLS model with a transition function $T(x)$ which smoothly, and monotonically varies between  zero and one.  Specifically, in this paper we take 
\begin{equation}
    \label{eq:transition}
    T\left( x \right)= \frac{1}{2} \left[ 1+ \frac{2}{\pi} \text{tan}^{-1}\left(\frac{\pi (x-c)}{w}\right) \right] ,
\end{equation}
where $w$ determines the steepness of the transition region. The approximate Heaviside function  used in the original PaLS work also monotonically increased from zero to one \cite{aghasi10} and thus is  similar to $T(x)$.  The difference lies in the width of the transition.  As seen in Figure \ref{fig:transition_a}, for the approximate Heaviside function, the transition region is by construction very narrow, as the goal was to represent binary valued objects, whose ``phases'' are separated by a narrow, 
smooth interface. In the new formulation, we stretch the width of this transition region, so that $f(\rmb)$ can assume basically \textit{any} value between $C_{L}$ and $C_{H}$, as seen in Figure \ref{fig:transition_b}. 

\begin{figure}[!ht]
\centering
\subfloat[]
{\includegraphics[width=0.32\textwidth]{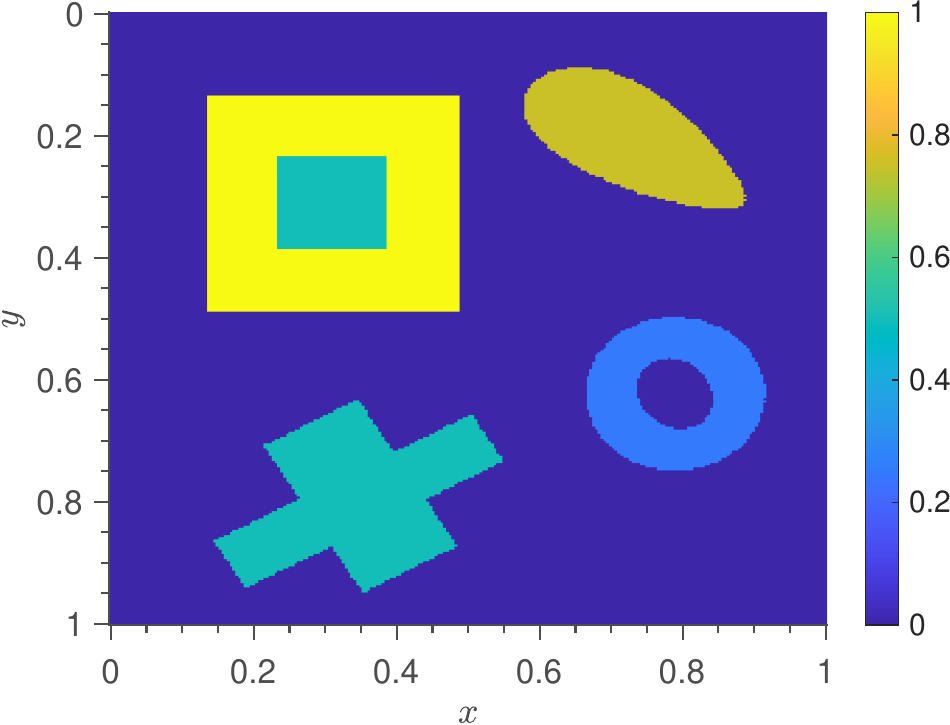}
\label{fig:mesh_a}}
\centering
\subfloat[]
{\includegraphics[width=0.32\textwidth]{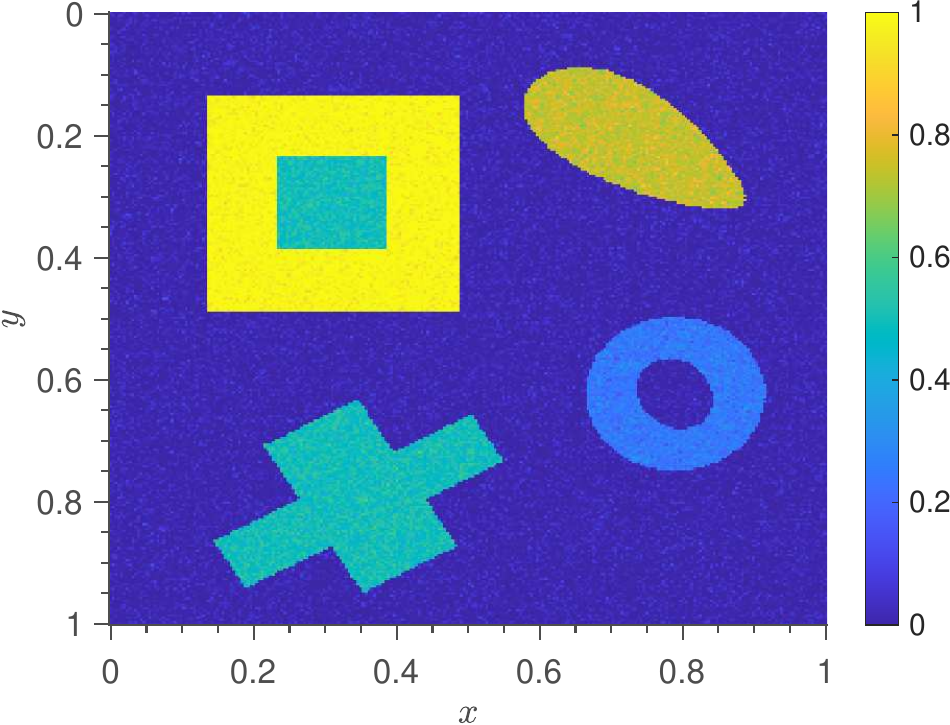}
\label{fig:mesh_b}}
\centering
\subfloat[]
{\includegraphics[width=0.32\textwidth]{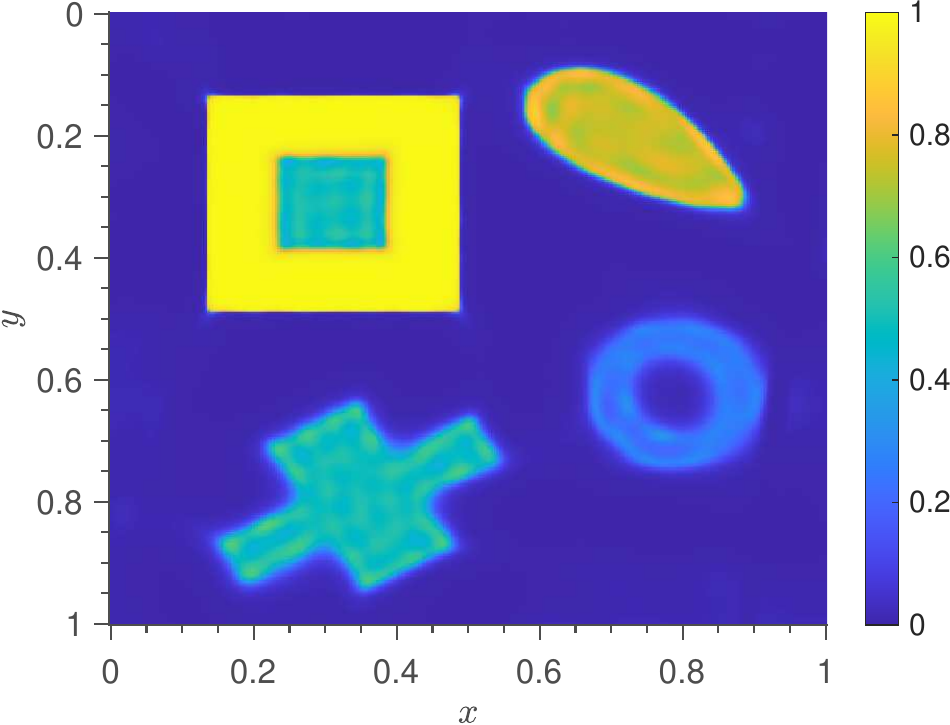}
\label{fig:mesh_c}}
\centering\\
\subfloat[]
{\includegraphics[width=0.32\textwidth]{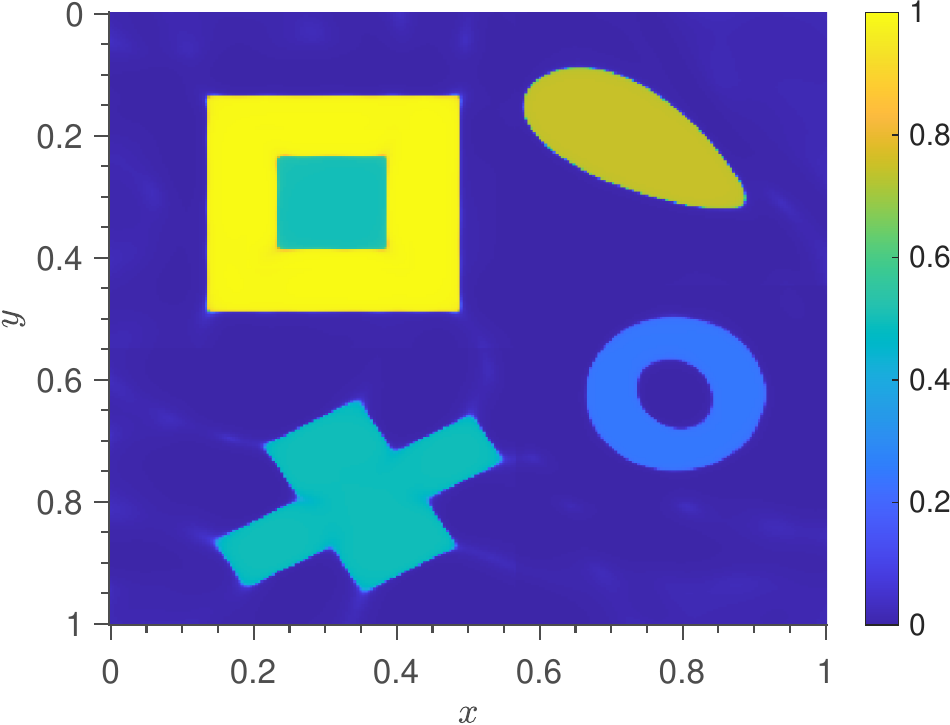}
\label{fig:mesh_d}}
\centering
\subfloat[]
{\includegraphics[width=0.32\textwidth]{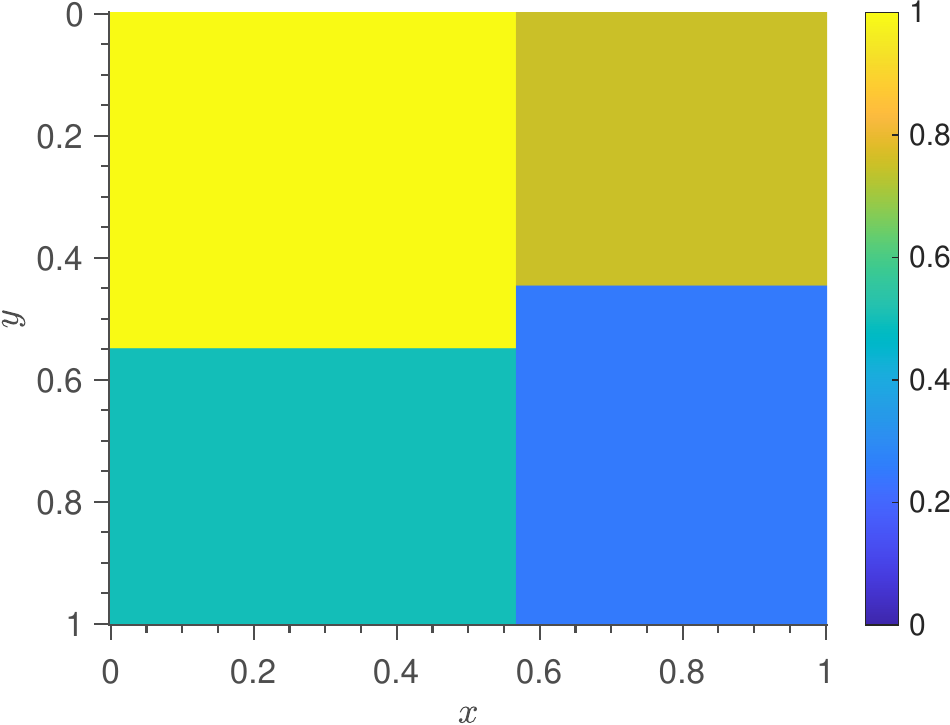}
\label{fig:mesh_e}}
\centering
\subfloat[]
{\includegraphics[width=0.32\textwidth]{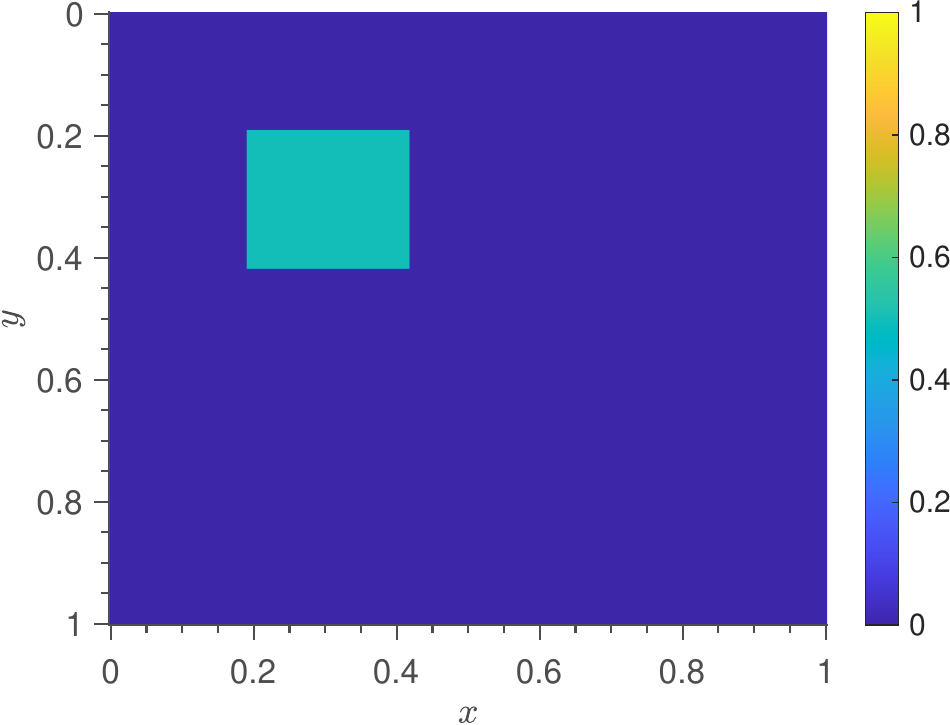}
\label{fig:mesh_f}}
\centering \\
\subfloat[]
{\includegraphics[width=0.32\textwidth]{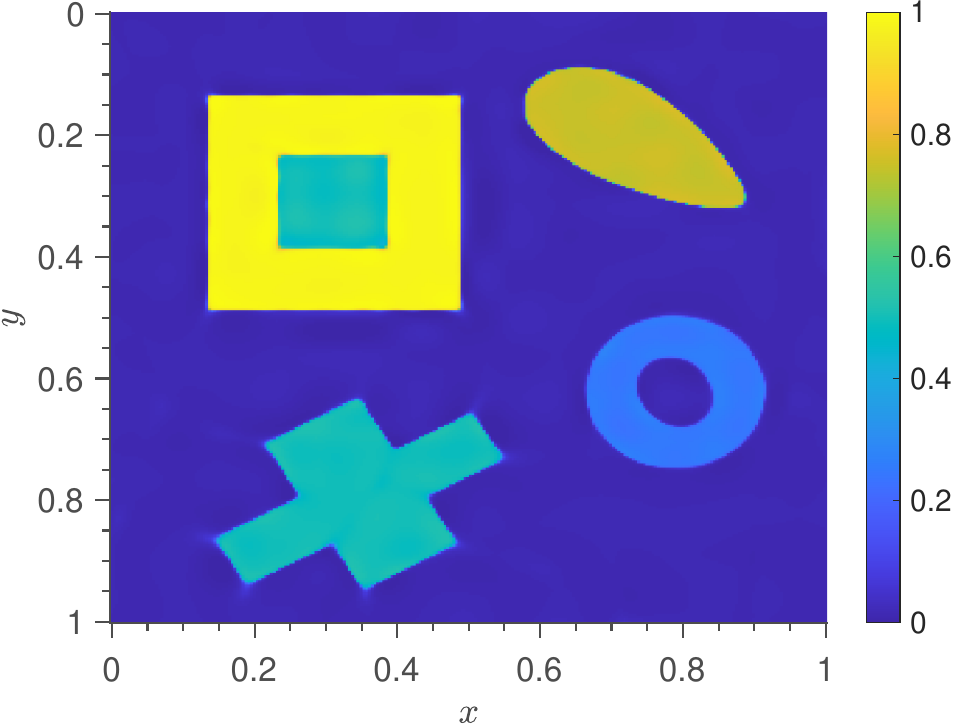}
\label{fig:mesh_g}}
\centering
\subfloat[]
{\includegraphics[width=0.32\textwidth]{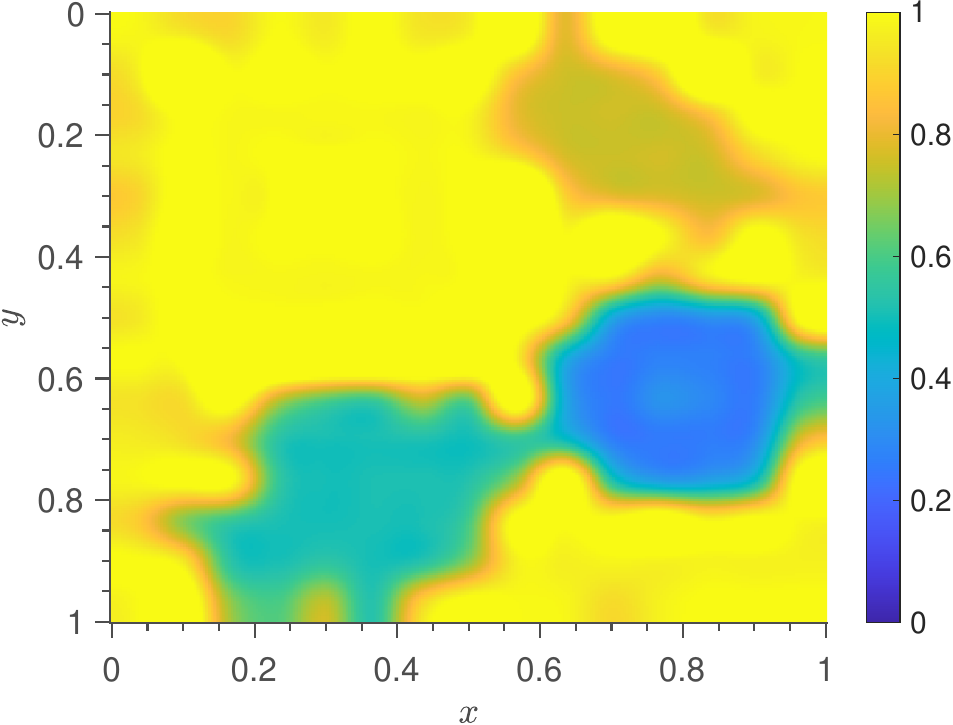}
\label{fig:mesh_h}}
\centering
\subfloat[]
{\includegraphics[width=0.32\textwidth]{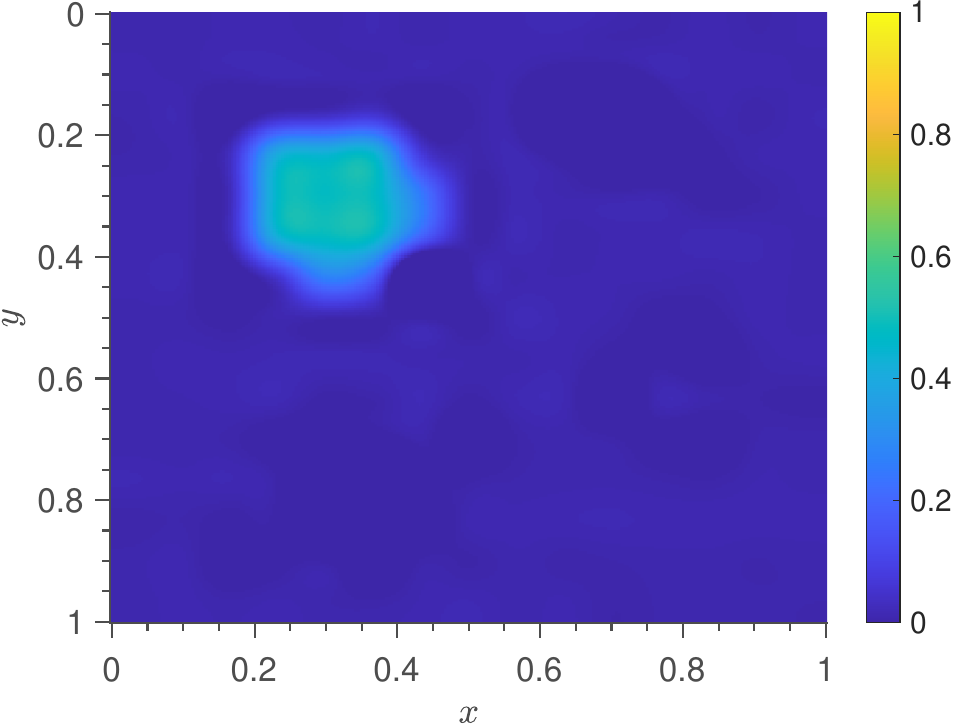}
\label{fig:mesh_i}}
\centering
\caption{(a) Image of size $256 \times 256$ pixels comprised of five piecewise constant objects with different shapes and contrasts on a zero-contrast background. (b) The image corrupted with additive Gaussian noise (SNR: 22dB). (c)  \palentir{} recovery of the image with fixed contrast limits. (d) \palentir{} recovery using manually adjusted vectors (e) $\Cmb_{H}$ and (f) $\Cmb_{L}$. (g) \palentir{} recovery with the parameterized contrast limits, $\mathbf{p}_c \in \R^{2N}$, included in the estimation problem. Estimated parameterized contrast limits (h) $C_{H}(\rmb)$ and (i) $C_{L}(\rmb)$.}
\label{fig:mesh}
\end{figure}

To illustrate the impact of these changes to the PaLS model, in Figure \ref{fig:mesh_a}, we display an image of size $256 \times 256$ pixels, comprised of five piecewise constant objects on a zero-contrast background. Here, we took the forward model to be the identity, so that $\dmb = \fmb(\pemb) + \wmb$, where $\fmb(\pemb)$, a discrete representation of $f(\rmb; \mathbf{p}, C_H, C_L)$, is a  vector of length $N_{pts} = 256^2 = 65536$ obtained using the discretization process described in section \ref{sec:experiment}. 
The vector $\wmb$ contains independent, identically distributed Gaussian random variables, with variance such that the signal to noise ratio (SNR) of the corrupted image is 22 dB. In this paper, SNR is calculated as
\begin{equation} \label{eq:snr} 
\text{SNR}_{\text{dB}}=20\text{log}_{10} \left( \frac{\sqrt{\sum^{N_{pts}}_{n=1}|d_n|^2}}{\sqrt{\sum^{N_{pts}}_{n=1}|w_n|^2}} \right)
\end{equation}
where $d_n$ and $w_n$ are the $n^{th}$ element of the noise-free data and the noise, respectively. The image corrupted with additive Gaussian noise is shown in Figure \ref{fig:mesh_b}, along with our \palentir{} reconstruction in Figure \ref{fig:mesh_c}. We used $225$ basis functions centered on an equally spaced $15 \times 15$ grid. Each basis function requires the parameters $\alpha$, $\beta$ and $\gamma$, resulting in total of $3N$ (675) parameters for $N$ (225) basis functions. Contrast coefficients are chosen as  $C_{H} = 1$ and $C_{L} = 0$,
since the maximum contrast in the image is equal to 1 and the minimum is equal to 0.

In Figure \ref{fig:mesh_c}, we see that the bright yellow region, whose contrast is equal to $C_{H}$, is recovered most accurately among the objects. The $\alpha_i$ of the basis functions in this region take on relatively large values which are truncated by $T(x)$. The regions without objects are also recovered very accurately. In those areas the coefficients are relatively small, so that the transition function $T(x)$ creates a constant, zero-background. For the remainder of the structures whose contrasts are between 0 and 4, we see blurred edges and oscillations similar to what we would expect from least squares denoising. Overall, the use of $T_{w}$ resulted in accurate recovery of regions with values not between $C_{H}$ and $C_{L}$, and rough recoveries of both shapes and contrasts of regions with values in between $C_{H}$ and $C_{L}$. This is achieved through the utilization of a single level-set, requiring only 675 unknowns, despite dealing with an 65536-pixel image.

The strong results in Figure \ref{fig:mesh_c} for the background and the yellow object in the upper left corner arise from the ability of the transition function to truncate the  values of $\phi $ which are not between $C_{H}$ and $C_{L}$. More generally, \palentir{} can achieve similar results across the entire scene by varying
${C}_{H}$ and ${C}_{L}$ in space.
For example, consider the case where we set $C_{H}(\rmb)$ and $C_{L}(\rmb)$
as shown in Figures \ref{fig:mesh_e} and \ref{fig:mesh_f}. Re-estimating $\pemb$ with these distributions of $C_{H}(\rmb)$ and $C_{L}(\rmb)$ yields the results in Figure \ref{fig:mesh_d}, where we now see near-perfect recovery of all objects and the background.

The results in Figure \ref{fig:mesh_d} and the associated discussion illustrate the ability of the \palentir{} model to recover piecewise constant scenes comprised of objects with more than two contrasts using a single level-set function.  This feature  rests heavily on the ability to specify space-varying bounds on the contrasts at relatively coarse scales. The key issue is to do this in a way that retains the advantages of a PaLS-type model: low-order and no need for explicit regularization. This is accomplished by parameterizing $C_{H}(\rmb)$ and $C_{L}(\rmb)$, the spatially varying contrast limits, using bi-cubic (tri-cubic for 3D) interpolation \cite{bicubic,tricubic}
(via the Matlab function ``\textit{imresize}”). In addition to $\alpha, \beta, \gamma$, for each basis function, we add 2 new parameters, which we refer as $\textit{p}_H$ and $\textit{p}_L$,  for upper and lower contrast bounds respectively. Consequently, for a \palentir{} model with $n_{B}\times n_{B}$ basis functions, where $n_B = \sqrt{N}$ , there are two $N$-vectors of parameters, 
referred as $\mathbf{p}_H$ and $\mathbf{p}_L$. 
To see the impact of parameterizing the contrast limits, we repeat the same experiment with the same $\wmb$ to recover the multi-contrast objects in Figure \ref{fig:mesh_a}. For the initial values of $\pemb$, we use the previously estimated \palentir{} parameters, which represent the image in Figure \ref{fig:mesh_c}, and we initialize $\mathbf{p}_H$ and $\mathbf{p}_L$ to ones and zeros respectively and rerun the \texttt{TREGS} algorithm. Figure \ref{fig:mesh_g} shows the resultant \palentir{} representation. The number of estimated parameters is increased from 675 to only 1125, still significantly fewer than the number of pixels in the represented image, which is 65536. As seen in Figure \ref{fig:mesh_g}, the new model with parameterized contrast limits has recovered all objects and the background accurately with near-perfect 
sharp boundaries, and maintains smooth piecewise constant contrasts everywhere. The contrast limits produced with the estimated $\mathbf{p}_c = 
\left[\mathbf{p}_H; \mathbf{p}_L \right]$, are shown in Figures \ref{fig:mesh_h} and \ref{fig:mesh_i}. The estimated contrast limits capture the contrast values of the objects at a relatively coarse scale. Noticeably, they look very similar to the handcrafted $C_{H}(\rmb)$ and $C_{L}(\rmb)$ in Figures \ref{fig:mesh_e} and \ref{fig:mesh_f}, especially at or near where the objects are located. 

\subsection{Advancements in the shape-expressiveness}
\label{sec:palentir_3}
While radial basis functions (RBFs) provide flexibility in terms of shape representation for PaLS, they are limited by the fact that they posses only circular, in 2D, or spherical, in 3D, level-sets, making them inefficient for representing, e.g., highly anisotropic shapes  \cite{est19}.  Motivated by these observations, we replace the RBFs in the old PaLS \cite{aghasi10} with a model of the form
\begin{equation}
    \psi (\rmb)=  e^{-\norm{\left( \Rmb (\rmb-\boldsymbol{\chi})\right)}^2_2}, 
    \label{eq:newpsi}
\end{equation}
where $\Rmb$ is a $2 \times 2$ matrix for 2D problems and $3 \times 3$ for 3D problems. As in \cite{est19}, this  model produces elliptical cross sections; however, the parameterization of $\Rmb$ is new and requires fewer parameters.  We discuss the 2D and 3D cases separately.  

In the 2D case, we define
\begin{equation} 
    \Rmb= \mu   \begin{bmatrix}
                    e^{\beta} & \gamma\\
                    0 & e^{-\beta} 
                \end{bmatrix} ,
    \label{eq:non_isotropic_rbf}
\end{equation}
where as shown below, $\beta$ and $\gamma$ define the eccentricity and orientation of the elliptical level-sets, while $\mu$ impacts the scale of the ellipses. In this paper, $\beta$ and $\gamma$ are parameters to be determined, while $\mu$ is fixed.  To elucidate the role of these parameters, we consider the $N=1$ case for both PaLS and \palentir{}.  We define the  $c$-level-set representation of a curve $\Gamma$ as $\Gamma= H(\phi(\rmb,\pemb)-c)$.  That is, $\Gamma$ is the set of points $\rmb$ such that $\phi(\rmb;\pemb)=c$ with $c > 0$: 
\begin{equation} 
    \Gamma = \left\{(x,y) \in \mathbb{R}^2 | \phi(x,y;\pemb)= c  \right\}.
    \label{eq:curve}
\end{equation}
We begin by examining the PaLS $c$-level-set of a single RBF centered at the origin:
\begin{equation}
    \phi_{\rm rbf} (x,y,\alpha,\beta)=\alpha  e^{-\beta(x^2+y^2)}.
    \label{eq:rbf}
\end{equation}
With $\phi_{\rm rbf}$ in \eqref{eq:curve} and $c>0$, we require $\alpha > 0$ as well to obtain a nonempty shape. Simple algebra gives
\begin{equation} 
    \Gamma_{\rm rbf} = \left\{(x,y) \in \mathbb{R}^2 \left| x^2+y^2=\frac{1}{\beta}\ln{\frac{\alpha}{c}} \right.    \right\}.
    \label{eq:curve_rbf}
\end{equation}
As anticipated, \eqref{eq:curve_rbf} defines a circle with radius $\sqrt{\frac{1}{\beta}\ln{\frac{\alpha}{c}}}$. We emphasize here  that, with $c$ fixed, there are an infinite number of $(\alpha, \beta)$ pairs that will give the same circle.

Next, we consider the $N=1$ \palentir{} $c$-level-set where 
$\phi(x,y,\alpha,\Rmb)=\sigma_h(\alpha)\exp\left(-\|\Rmb\rmb\|^2_2\right)$ with $\rmb^T = [x\;y]$.  In this paper, we take $\sigma_h(\alpha) = \tanh{\alpha}$.
Notice that $\sigma_h(\cdot)$ bounds the weight coefficient of the basis function between $-1$ and $1$ and is discussed in greater depth in section \ref{sec:palentir_4}. Using this $\phi$ in \eqref{eq:curve} and $c>0$, a non-empty  $\Gamma$ requires $\sigma_h(\alpha)>0$ and thus $\alpha > 0$.  Dividing both sides by $\sigma_h(\alpha)$ and taking the logarithm yields
\begin{equation}
    \Gamma_{\rm new} =\left\{(x,y) \in \mathbb{R}^2 \; |\;  \left(e^{\beta}x+\gamma y \right)^2 + e^{-2\beta}y^2=\tau^2  \right \}
    \label{eq:curve_rbf_new}
\end{equation} 
with $\tau^2 = \frac{1}{\mu^2}\ln{\frac{\sigma_h(\alpha)}{c}}$.  When $\beta$ and $\gamma$ are $0$, similar to $\Gamma_{\rm rbf}$, $\Gamma_{\rm new}$ is a circle with a radius of $\tau$.  As $0 < \sigma_h(\alpha) \leq 1$, we see that $\mu$ and $c$  define the largest circle representable using \palentir{}.  While 
it may be useful to estimate
one or both of these quantities along with the other model parameters, here we choose them 
to be constants.  Specifically, $c$ is set to $0.01$, a value slightly larger than $0$, for the reasons discussed in \cite{aghasi10}, and $\mu$ is set to $10$ to have the maximum area of a single ABF to be slightly less than $15\%$ of the total area, to restrict
the extent to which the ABFs overlap each other.

\begin{figure}[!ht]
    \centering
    \subfloat[]
        {\includegraphics[width=0.49\textwidth]{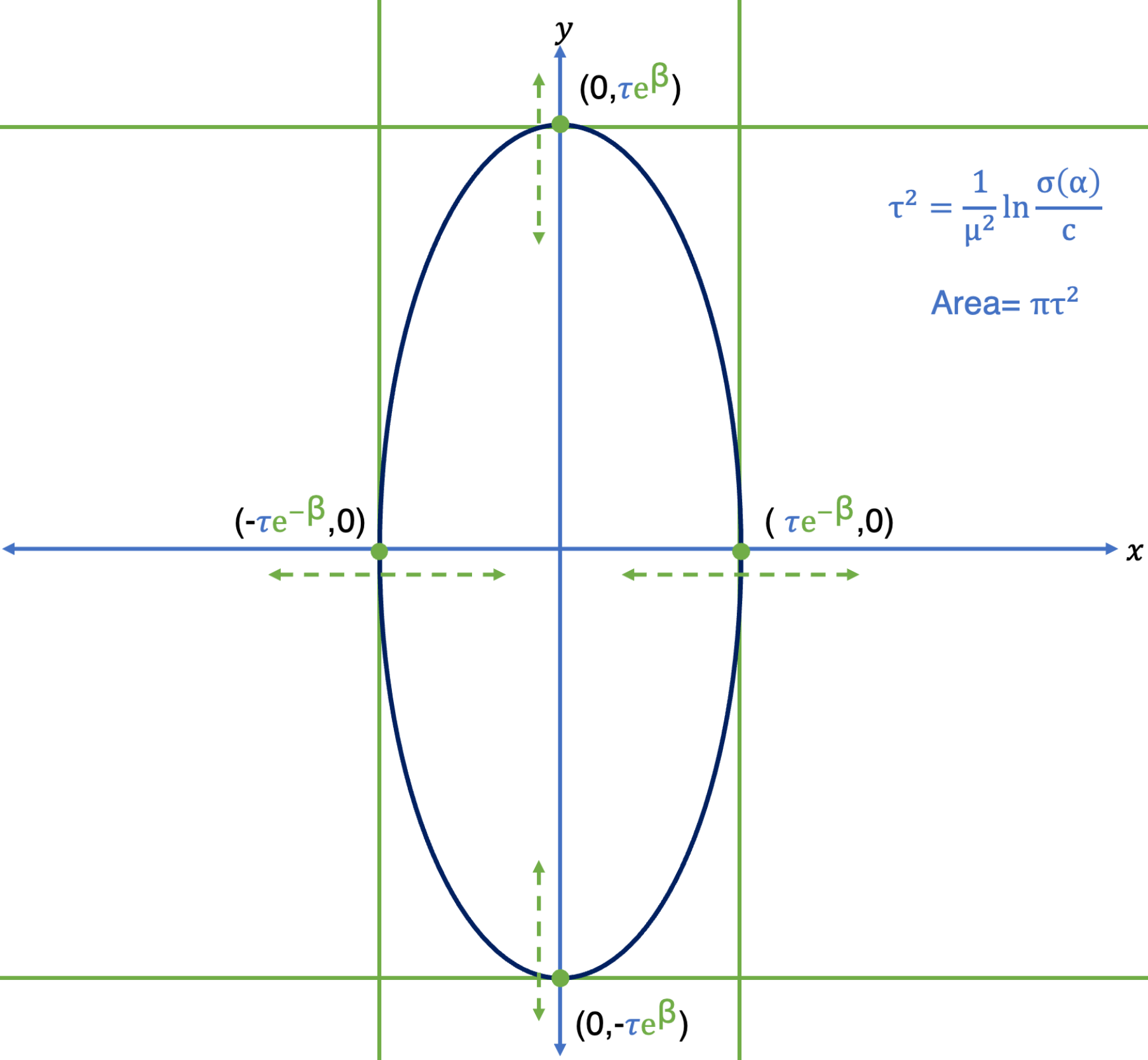}
        \label{fig:clevel_ellipse_a}}
        \centering
        \subfloat[]
         {\includegraphics[width=0.49\textwidth]{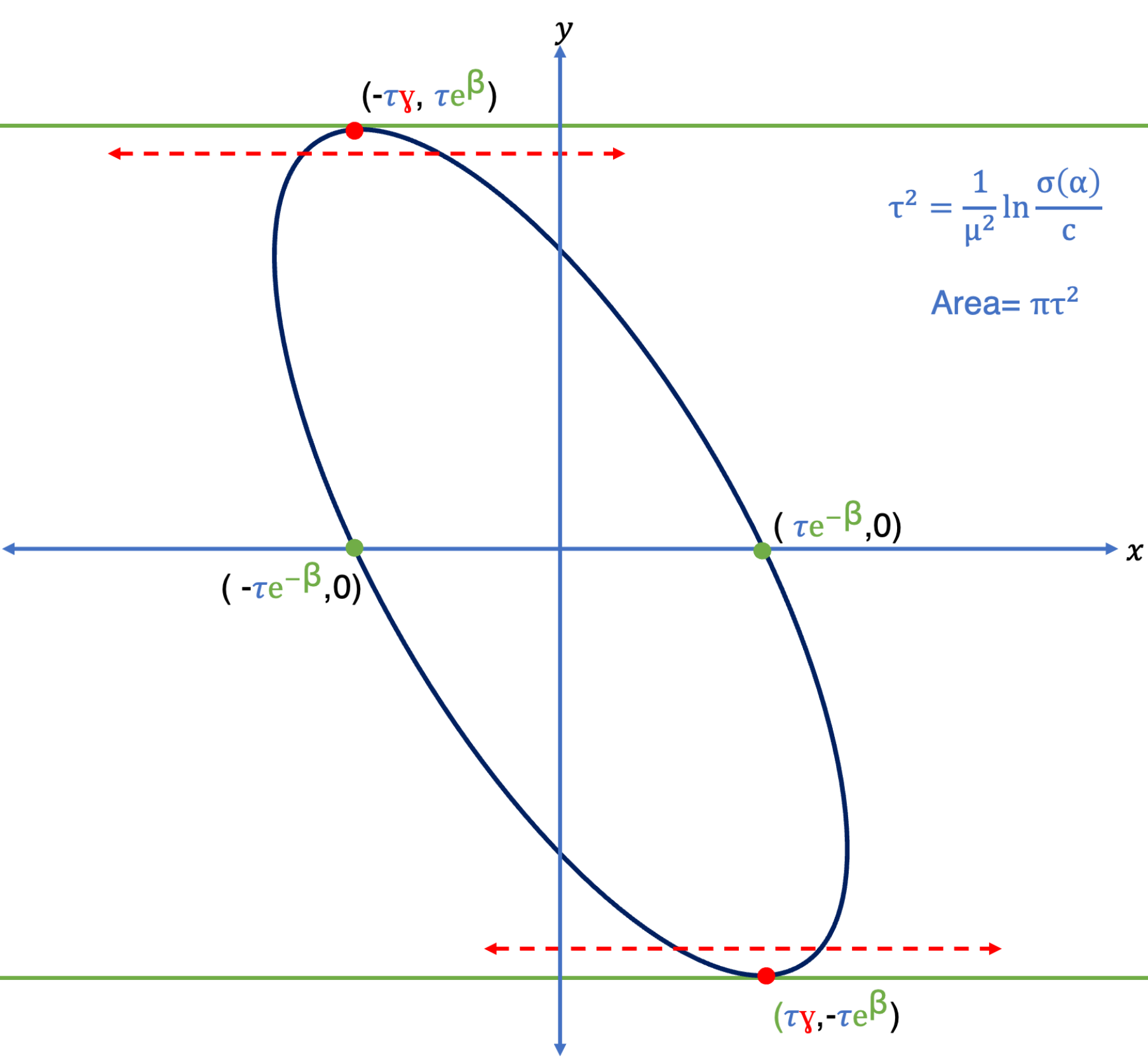}
        \label{fig:clevel_ellipse_b}}
                 \centering
    \caption{The impact of the parameters on the $c$-level ellipsoidal representation of the new basis function. (a) shows the impact of the \textit{stretching} parameter $\beta$ and (b) shows the impact of the \textit{sliding} parameter $\gamma$.}
         \label{fig:clevel_ellipse}
\end{figure}

When $\beta$ and $\gamma$ are nonzero, the $\Gamma_{new}$ curve becomes elliptical. To understand the role of each parameter in defining the geometry of the ellipse, we start by viewing $\Rmb$ as the Cholesky factor of the symmetric positive definite (SPD) matrix $\Amb = \Rmb^T\Rmb$, which we call the \textit{Stretch and Slide Matrix}. This matrix includes two types of parameters which are estimated by the \palentir{} reconstruction: the \textit{stretching} parameter $\beta$ and the \textit{sliding} parameter $\gamma$.  In Figure \ref{fig:clevel_ellipse_a}, the impact of the \textit{stretching} parameter $\beta$ is shown. With $\gamma = 0$, nonzero values of $\beta$ produce ellipses with principle axes in the $x$ and $y$ directions. We visualize the impact of increasing $\beta$ as someone holding the circle (assuming $\beta = 0 $ initially) at two opposite points touching the $x$ and $y$ axes (points shown with green dots in Figure \ref{fig:clevel_ellipse_a}), and stretching in the direction of green dashed arrows.
We call this shape transformation  ``stretching,'' and hence $\beta$ the \textit{stretching} parameter. Similarly, as shown in Figure \ref{fig:clevel_ellipse_b}, $\gamma$ is the \textit{sliding} parameter, because increasing $\gamma$ horizontally ``slides'' the maximum and the minimum points of the ellipse along the vertical axes, as shown in  with red dashed arrows. 
Since the determinant of $\Amb$ is $\mu^2$, these two shape transformations do not change the area of the ellipse. The area depends on the constants $\mu$, $c$, and the parameter $\alpha$. In summary, with $c$ and $\mu$ fixed, $\alpha$ controls the total area by homogeneously expanding the shape, $\beta$ controls the upper and lower tangent lines by stretching the ellipse, and $\gamma$ slides the tangent points of the ellipse along the tangent lines.

If we compare \eqref{eq:curve_rbf} and \eqref{eq:curve_rbf_new}, we can clearly see that the $\beta$ in $\phi_{\rm rbf}$, which adjusts the scale of the circles, is replaced with $\mu^2$. As we explained above, the parameters $\alpha$ and $\beta$ in $\phi_{\rm rbf}$ have a related (but opposite) effect on the $c$-level-set, and infinitely many pairs ($\alpha, \beta$) exist 
that give the same $\Gamma_{\rm rbf}$. Hence, we replaced the $\beta$ in $\phi_{\rm rbf}$ with a fixed constant $\mu^2$ to reduce the dimension of the search space of the model and use only the parameter $\alpha$ to determine the scale of the ellipses.

\begin{figure}[!ht]
     \centering
     \subfloat[]
        {\includegraphics[width=0.32\textwidth]{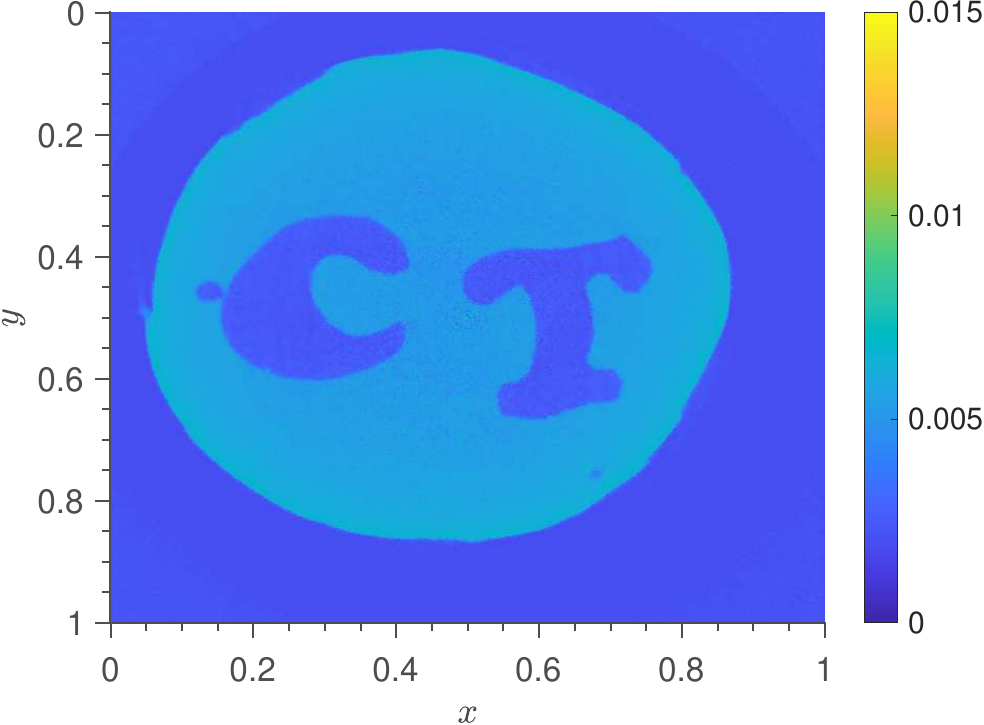}
        \label{fig:cheese_GT1}}
     \centering
     \subfloat[]
         {\includegraphics[width=0.32\textwidth]{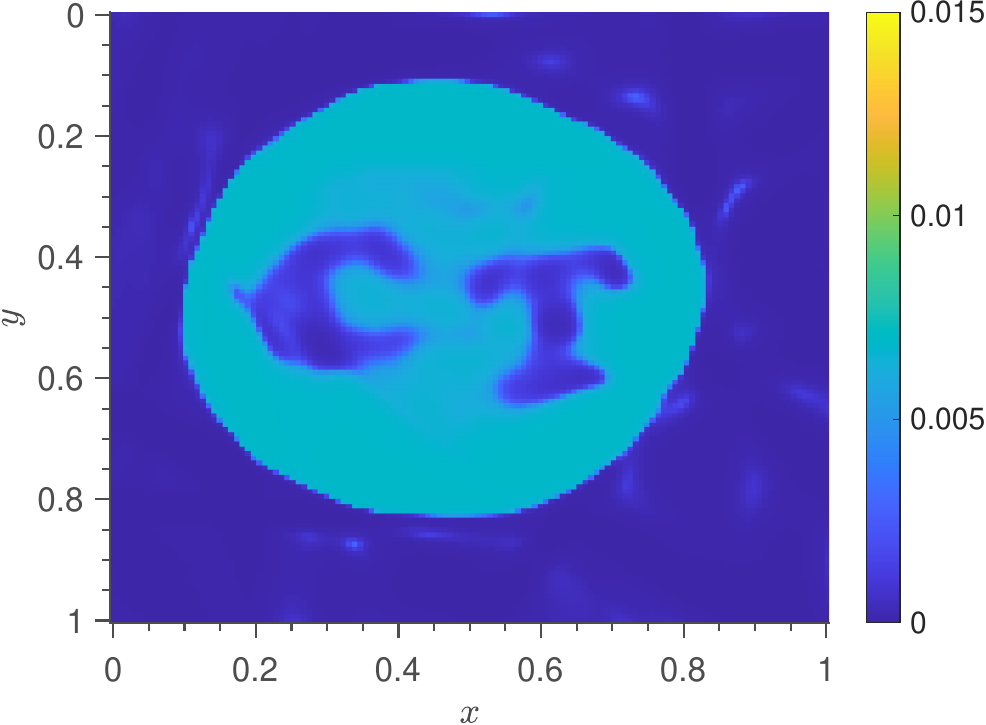}\label{fig:cond_cheese_palentir}}
         \centering
         \subfloat[]
         {\includegraphics[width=0.32\textwidth]{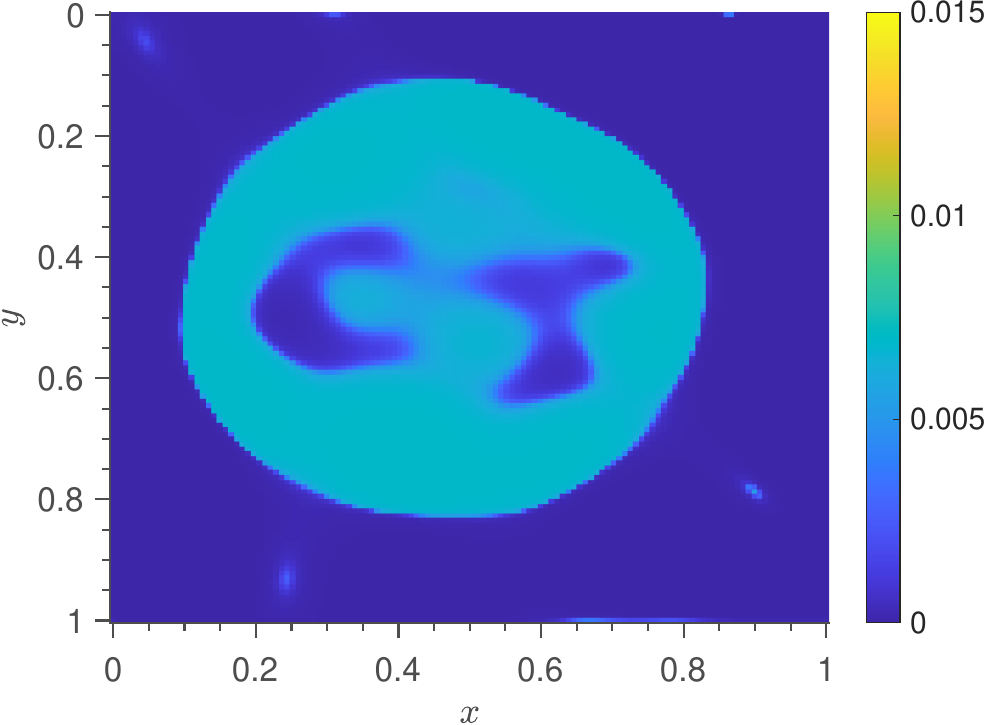}\label{fig:cond_cheese_pals}}
         \centering
    \caption{(a)  Reference image of 
    a carved cheese specimen
    computed  using a high-resolution filtered back-projection (FBP) 
    from the 360-projection sinogram. (b) \palentir{} recovery from the tomographic X-ray data using only 15 projections.  (c) RBF PaLS recovery 
    using the same data as \palentir{}.
    }
    \label{fig:cheese_cond}
\end{figure}

To demonstrate the extended capacity for shape representation 
in \palentir{}, we compare
\palentir{} and RBF PaLS on a 2D X-ray computed tomography (CT) test problem.
Utilizing an open-access data set of tomographic X-ray data, featuring a carved cheese specimen \cite{bubba_cheese}, the reference image produced through high-resolution filtered back-projection (FBP) reconstruction, yielding a $2000\times 2000$ pixel representation, is illustrated in Figure \ref{fig:cheese_GT1}. Focusing on a sparse view problem with 15 projections spanning the  full 360 degree circle, both \palentir{} and RBF PaLS configurations used $8\times 8$ basis functions. Additionally, both models featured fixed values for the contrast coefficients, $C_{H}$ and $C_{L}$ for \palentir{} and $f_O $ and $f_B $ for RBF PaLS, set to $0.007$ and $0$, respectively. The recovered images of $128\times 128$ pixels, are shown in Figures \ref{fig:cond_cheese_palentir} and \ref{fig:cond_cheese_pals}. As discussed earlier, 
using the new ABFs introduces greater adaptability in shape representation, allowing for the capture of finer details with the same number of basis functions, and in fact reducing the total number of unknowns as well. In 2D, for $N$ basis functions, RBF PaLS requires $4N$ parameters, $\alpha_j$, $\beta_j$, and $\boldsymbol{\chi}_j$, for $j=1,2,...N$. Whereas, for $N$ basis functions \palentir{} requires only $3N$ parameters, $\alpha_j$, $\beta_j$, and $\gamma_j$. The results underscore the advantages of \palentir, 
as it exhibits a more detailed recovery, discerning subtle features that 
are not captured
in the RBF PaLS representation. For instance, the small hole on the left side of the carved letter ``C'' is effectively captured in the \palentir{} recovery, albeit not perfectly, while the hole cannot be observed 
in the RBF PaLS recovery. Similarly, when analyzing the carved letters ``C'' and ``T,'' the curved edges and fine details are more faithfully restored in the \palentir{} representation, while they are smeared out
in the RBF PaLS recovery. 
Overall, RBF PaLS was able to recover the main body of the cheese and the background, 
but the recovery of the carved letters is rough. 
On the other hand, \textit{with the same number of basis functions and 
fewer (64 less) 
parameters}, \palentir{} recovers 
significantly more details on the carved letters and even the small details on the cheese, which are 
absent from the RBF PaLS recovery. 

For the 3D case, we construct $\Rmb_3 \in \R^{3\times 3}$ using three stretch and slide matrices as follows
\begin{align*}
\label{eq:non_isotropic_rbf_3d}
   &\Rmb_{3}=\\&  \mu \begin{bmatrix}
        e^{\beta_{1}} & \gamma_{1} & 0 \\
        0 & e^{-\beta_{1}} & 0 \\
        0 & 0 & 1
    \end{bmatrix} \times
    \begin{bmatrix}
        1 & 0 & 0 \\
        0 & e^{\beta_{2}} & \gamma_{2} \\
        0 & 0 & e^{-\beta_{2}}
    \end{bmatrix} \times
    \begin{bmatrix}
        e^{\beta_{3}} & 0 & \gamma_{3} \\
        0 & 1 &   0 \\
        0 & 0 & e^{-\beta_{3}}
    \end{bmatrix}.
\end{align*}
Note there are six parameters in our model, with each of the three matrices in the product affecting a shear transformation in a 2D plane within $\Rbb^3$. Similar to 2D case, $\Rmb_3$ creates anisotropy in $x,y,z$ and, assuming all three $\beta_i$ are finite, may be viewed as the Cholesky factor of the 
symmetric positive definite (SPD) matrix $\Amb = \Rmb_3^T\Rmb_3$. 
We note that in the 3D PaLS model in \cite{est19}, the authors 
modify the CSRBF formulation to use a matrix-based dilation to define ellipsoids.  Unlike what we propose here, their method uses a symmetric $3 \times 3$ matrix, parameterized by the six unique elements of such a structure.
One of the advantages of our formulation is that $\Rmb_{3}^T\Rmb_3$ is always SPD, whereas in the 3D PaLS model in \cite{est19} regularization is required to enforce this constraint.

\subsection{Numerical improvements and stability}
\label{sec:palentir_4}

As we detail below, relative to RBF PaLS, \palentir{} exhibits enhanced numerical performance through several key modifications: (a) a reduction in the number of parameters achieved by fixing the centers of the Anisotropic Basis Functions (ABFs); (b) the introduction of the $\sigma_h(\cdot)$ function to ensure that the values of $\alpha_i$ are bounded; and (c) alterations to the basis functions, ensuring uniqueness in the parametric representation of shapes by a better choice
of basis functions. 

In contrast to RBF PaLS, which requires adjusting RBF centers due to the limited shape expressiveness of RBFs, \palentir{} exhibits enhanced shape expressiveness with anisotropic basis functions (ABFs) using fixed and homogeneously distributed centers. This not only eliminates the need for center estimation during inversion but also reduces the dimension of the search space, contributing to numerical improvements.

Here, we use the condition number of the Jacobian matrix to
quantify numerical performance.  Recall that the condition number measures the sensitivity to perturbations of the solution to a linear system of equations. As discussed in \cite{jorge06, higham2002}, this is a crucial metric in determining the performance of a Newton-type approach to the estimation of the model parameters. We define the condition number of the Jacobian as the ratio of the largest to smallest singular value. Values of the condition number near 1 indicate a well-conditioned matrix, and large values indicate an ill-conditioned matrix. The expressions detailing the derivatives of the \palentir{} parameterization with respect to model parameters, which are necessary for the derivation of the Jacobian, are provided in the \ref{sec:derivatives}. Note that we are solving a least squares problem, so the condition number of the least squares problem is not just the ratio of singular values of the matrix but the residual.

\begin{figure}
     \centering
     \subfloat[]
      {\includegraphics[width=0.32\textwidth]{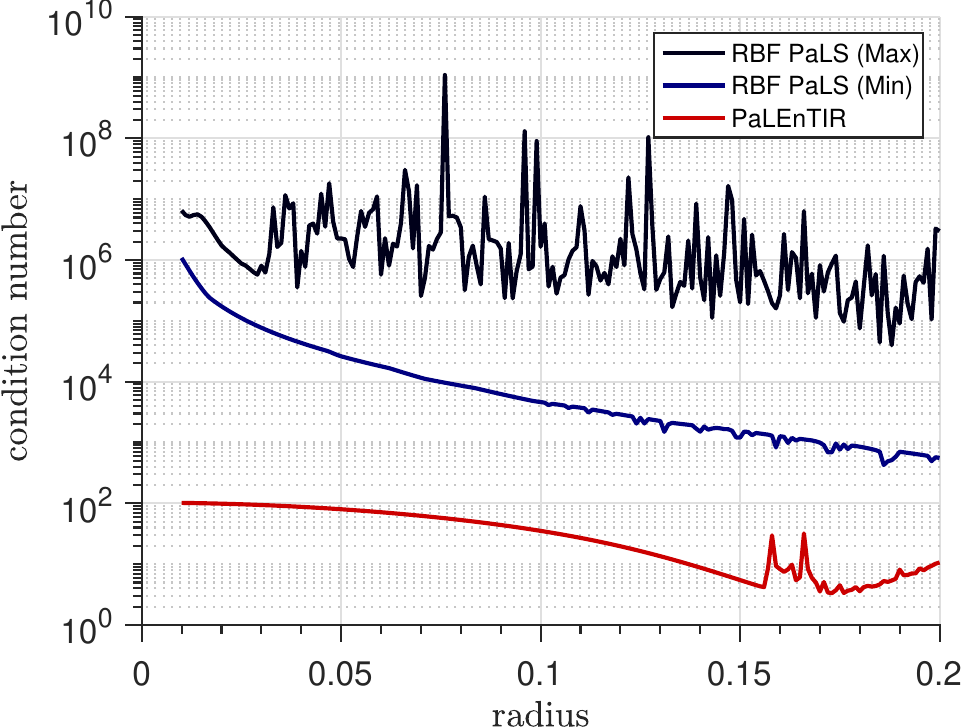}\label{fig:jac_a}}
     \centering
     \subfloat[]
      {\includegraphics[width=0.32\textwidth]{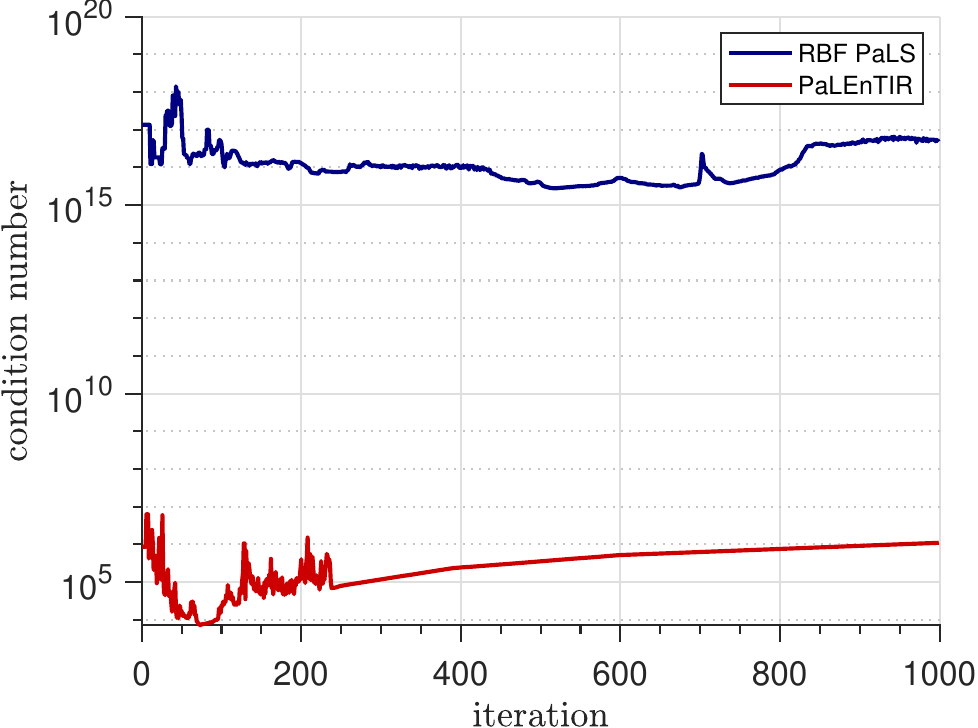}\label{fig:jac_b}}
         \centering
         \subfloat[]
         {\includegraphics[width=0.32\textwidth]{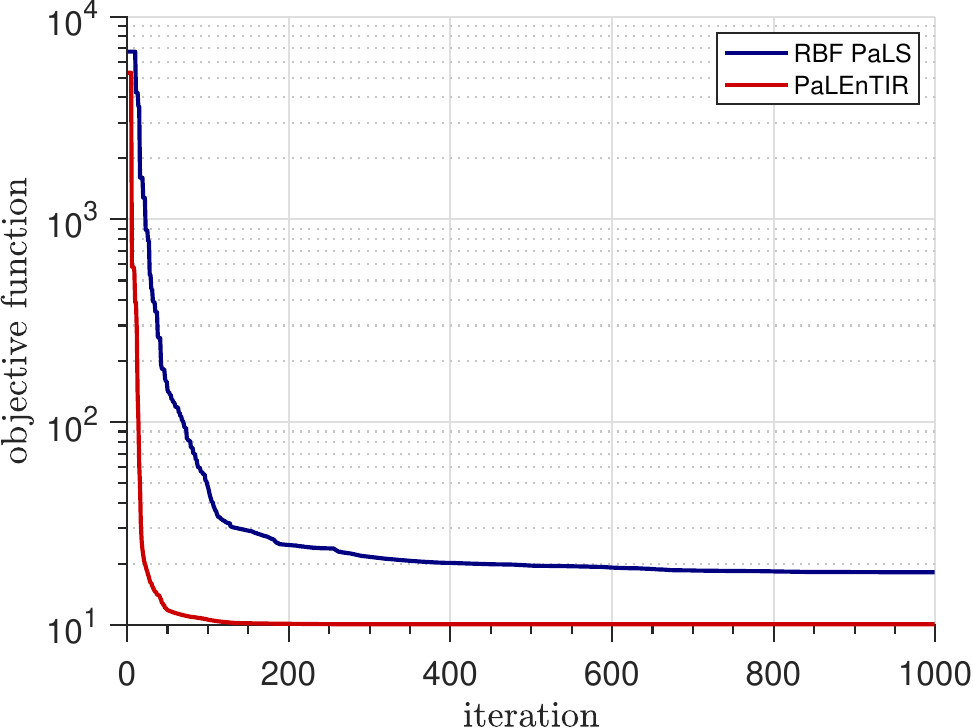}\label{fig:jac_c}}
         \centering
    \caption{(a) Condition number of the residual versus the radius of circular cross sections (corresponding to appropriate level sets) for the two models. The black line (Max) and blue line (Min) represent RBF PaLS; the  red line represents \palentir{}. (b) Condition number of the residual versus  \texttt{TREGS} iterations, and (c) objective function plots over  \texttt{TREGS} iterations of the 2D sparse view tomographic X-ray experiment. Blue line is for RBF PaLS and red line is for the new \palentir{}. The resultant reconstructions are shown in Figures \ref{fig:cond_cheese_palentir} and \ref{fig:cond_cheese_pals}.}
    \label{fig:jac}
\end{figure}

Similar to the discussion in subsection \ref{sec:palentir_3}, our initial focus is on the single-basis function scenario. Given that RBF PaLS is confined to generating circular cross-sections (see section \ref{sec:palentir_3}), we evaluate the numerical performance of both models in terms of their condition numbers, when their respective basis functions shape identical circles in their $c$-level-set representations, and the forward model is taken to be the identity in a noise-free case. Figure \ref{fig:jac_a} plots the condition number against the radius of the circular cross-sections for both RBF PaLS and \palentir{} models. Given that the leading coefficient of \palentir{} limits the size of produced shapes, we focus on circular cross-sections with radii ranging from 0.01 to 0.2 units within a square region of 2 units per side. The condition number is computed using the function "\textit{cond(.)}" in Matlab. Unlike \palentir{}, RBF PaLS can produce the same cross-section with many $(\alpha, \beta)$ pairs. Consequently, we produce identical cross-sections for radii ranging from 0.01 to 0.2, with increments of 0.001. We vary $(\alpha, \beta)$ pairs, increasing $\alpha$ by 0.01 within the range of 0.01 to 1. For each radius, we plot both the minimum and maximum condition numbers. The black and blue lines represent the maximum and minimum condition numbers for RBF PaLS, respectively, and the red line represents the condition number for 
\palentir{}. Notably, even if specifically selecting the $(\alpha, \beta)$ pairs that yield the minimum condition number for RBF PaLS, this is still much larger than the condition number obtained with 
\palentir{}.
In the \textit{worst} case scenario for RBF PaLS, where we pick the $(\alpha, \beta)$ pairs with the largest condition number 
for each cross section radius, 
the condition number is significantly larger
than for \palentir{}.

Next, we turn our attention to evaluating the numerical performance for 
a 2D X-ray computed tomography experiment, as detailed in subsection \ref{sec:palentir_3}. Figures \ref{fig:jac_b} and \ref{fig:jac_c} plot 
the condition number and the objective function (or misfit) against the number of \texttt{TREGS} iterations, respectively. In both plots, the black line represents RBF PaLS, while the blue line represents \palentir{}. Figure \ref{fig:jac_b} clearly shows the superior conditioning achieved by \palentir{}, confirming the results observed in the single-basis function case. This analysis is conducted with an $8 \times 8$ grid of basis functions for a real sparse-view X-ray CT problem. Additionally, Figure \ref{fig:jac_c} highlights two noteworthy enhancements in \palentir{}. Firstly, \palentir{} attains significantly better objective function values, supporting the observations drawn from Figures \ref{fig:cond_cheese_palentir} and \ref{fig:cond_cheese_pals} in the previous subsection. Secondly, the \palentir{} parameterization leads to better solutions (smaller objective function values) in substantially fewer iterations 
compared with the RBF PaLS
parameterization, leading to a substantial improvement in numerical performance as well.

\section{Experimental evaluation}
\label{sec:experiment}

We explore the utility of \palentir{} using a variety of linear and nonlinear forward models, particularly for problems where the data are limited.  
All experiments require that we use a discrete form of the forward model developed in section \ref{sec:formulation}.   To keep the discussion simple, 
we use quadrature rules 
for all integral operators 
and finite differences 
if the forward operator involves a PDE.  We could use other methods of discretization, but these assumptions permit a straightforward extension of the continuous to discrete notation for purposes of this paper; that is, the vector of unknowns represents the values of the desired function $f(\bfr; \pemb)$ at a finite set of grid points.  Specifically, 
let ${\bfr}_i,$ for $i=1, \ldots,N_{pts}$, denote a discrete set of 2D or 3D spatial grid points in $\Omega$. 
We define the 
$N_{pts}$-vector $\bfpp$ as $[\bfpp]_i = f(\rmb_i;\bfp)$. In this discrete case, $ {\mathcal M} (\bfpp)$  represents the measured data for all sources given the discrete values $f(\rmb_i;\bfp)$.  The discrete forward model thus becomes
\begin{equation} 
\label{eq:discreteforward} 
    \bfd = \mathcal{M}(\bfpp) + {\bf w}.
\end{equation}

We compare the performance of the \palentir{} model with that obtained using the L2-Total Variation (TV) regularization. For the TV method we use \texttt{TVReg} in \cite{hansen}. As a distinct advantage over pixel-based methods, \palentir{} eliminates the need for a regularization parameter, in contrast to the TV method, where selecting a regularization parameter is essential. {\em In this paper, we always choose the TV regularization parameter that results in the minimal Mean Squared Error (MSE) when comparing with \palentir{}. This ensures that the TV result used in the comparison represents the best-case scenario, which generally would not be achievable in practice.} The number of unknowns employed by \palentir{} in our experiments is quite modest in comparison with the necessary number of unknowns for the TV method. For $N$ ABFs, \palentir{} employs $3N$ (for 2D) or $7N$ (for 3D) PaLS parameters, and additionally $2N$ contrast parameters. Consequently, the parameter vector $\pemb$ of \palentir{} is of size $5N$ (for 2D) or $9N$ (for 3D). On the other hand, L2-TV regularization is a pixel-based method and the number of unknowns is equivalent to the number of pixels in the discretized image.  Note that number of parameters for \palentir{} is independent of the resolution and depends only on the number of basis functions, one of the key advantages of PaLS-type methods.

We show results for both \palentir{} and TV methods and compare them based on the following metrics: PSNR (dB), SNR (dB), SSIM, and MSE. The Structural Similarity Index (SSIM) was created to replicate the human visual perception system \cite{SSIM}.  A value closer to 1 indicates that the two images are very similar, whereas a value closer to -1 indicates the opposite. Mean Squared Error (MSE) is a commonly used metric to quantify the difference in the values of corresponding pixels between the sample and the reference images. PSNR is a commonly used metric to quantify the signal quality by comparing the peak level of a desired signal to the level of noise. 

\subsection{Deconvolution}
\label{sec:experiment_6}

Deconvolution is a linear inverse problem that recovers a desired signal from its convolution  
with a filter or a distortion function associated with an instrument or the physics of the problem.
The discrete forward model is defined as in \eqref{eq:discreteforward} 
with $\bfpp$ as in 
\eqref{eq:newpals}, giving the   forward mapping $\mathcal{M}(\bfpp) = \Amb\bfpp$. The matrix $\Amb$ comes from discretizing an integral equation that represents convolution.

\begin{figure}[!ht]
     \centering
     \subfloat[]
        {\includegraphics[width=0.32\textwidth]{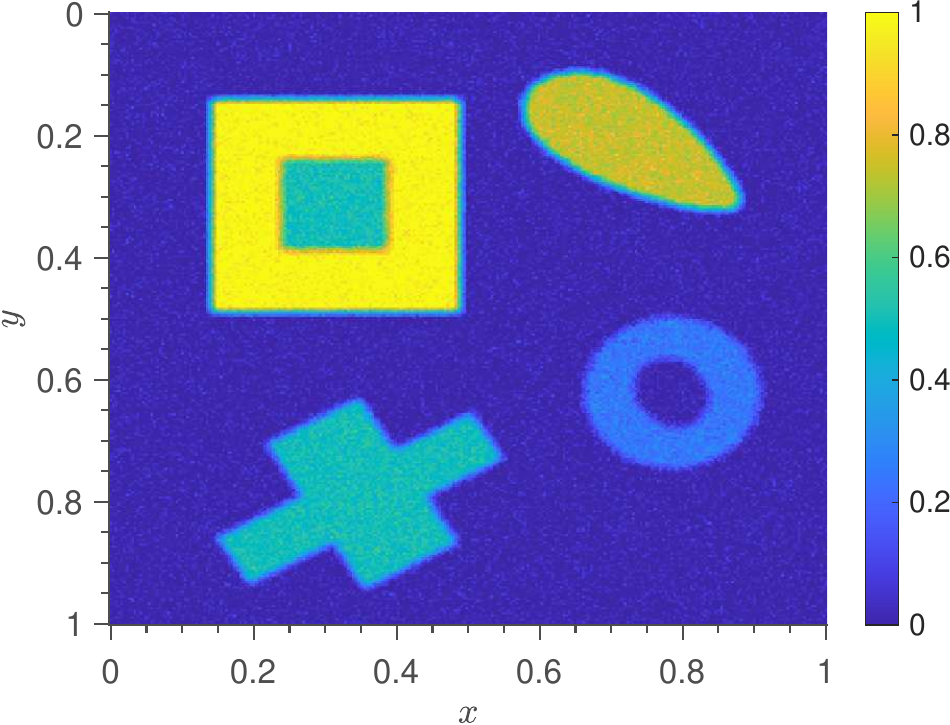}
        \label{fig:deconv_blur}}
                 \centering
                 \subfloat[]
        {\includegraphics[width=0.32\textwidth]{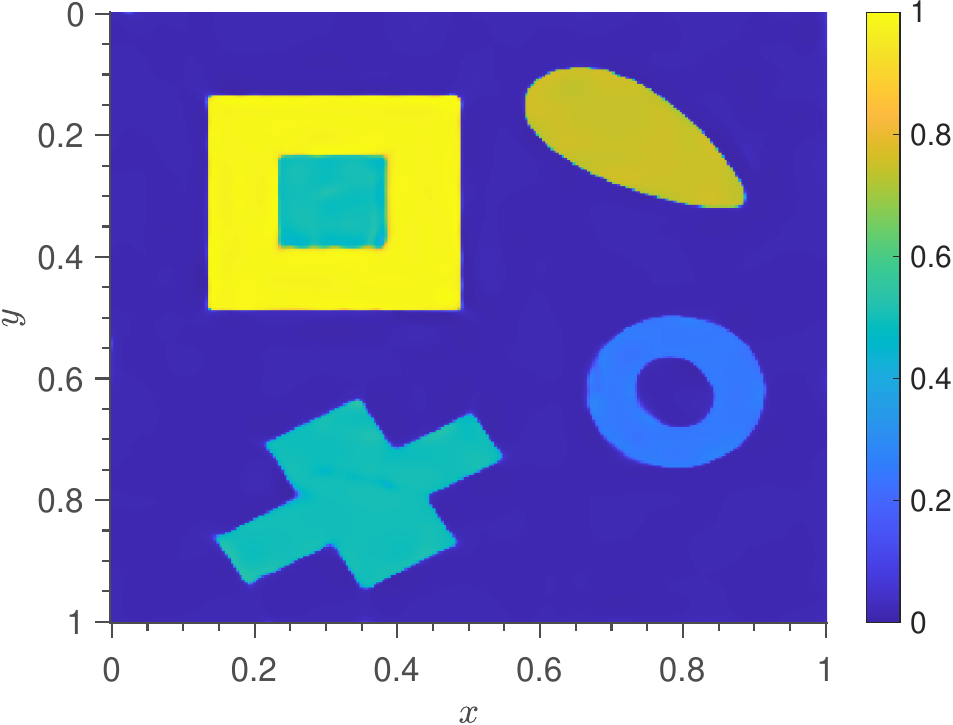}
        \label{fig:deconv_palentir}}
                 \centering
                 \subfloat[]
        {\includegraphics[width=0.32\textwidth]{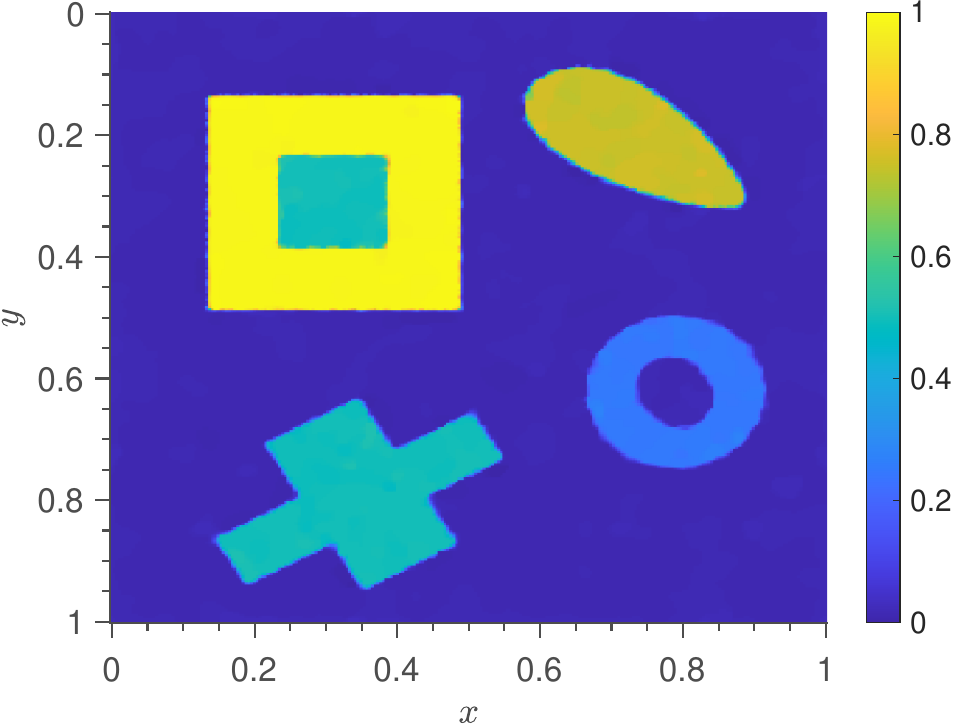}
        \label{fig:tv_deconv}}
                 \centering
\caption{(a) Blurred and noisy image of four objects. (b) \palentir{} reconstruction and (c) TV reconstruction for the deconvolution problem.}
\label{fig:deconvolution}
\end{figure}

In our deconvolution experiment, we seek to reconstruct a $256 \times 256$ image from input data that has been both filtered and corrupted by noise. The data is generated by convolving the authentic image with a rotationally symmetric Gaussian low-pass filter  of size $5 \times 5$ with standard deviation $1$ , using the MATLAB function \texttt{fspecial} with default parameters. Subsequently, the data is contaminated with $10 \%$ additive Gaussian white noise, resulting in a data SNR of $22$dB, as illustrated in Figure \ref{fig:deconv_blur}. We use 225 basis functions 
in a $15 \times 15$ grid resulting in 1125 unknown parameters.  This
results in a reduction of more
than 98\% in the number of unknowns compared with the $256 \times 256$ array of pixels comprising the underlying discretization of the image domain. 
Figure \ref{fig:deconv_palentir} shows the \palentir{} recovery, and Figure \ref{fig:tv_deconv} shows the TV recovery. Table \ref{table:deconvolution} shows the performance metrics of both recoveries; the best results are shown in bold font. \palentir{} performs better across all 4 metrics, while requiring relatively few unknowns. The larger PSNR and SNR and smaller MSE values indicate that \palentir{} outperforms TV in terms of minimizing the absolute error. Furthermore, according to the SSIM scores, the perceptual difference between the original image and the TV recovery is larger than that for \palentir{}, which suggests that compared with TV, the \palentir{} recovery exhibits a closer resemblance to the ground truth image. Upon reviewing Figures \ref{fig:deconv_palentir} and \ref{fig:tv_deconv}, it is evident that \palentir{} demonstrates a strong performance in delineating object boundaries, while the TV recovery exhibits some blurriness along several boundaries. Additionally, a close inspection of the TV recovery of the ring-shaped object reveals challenges in accurately representing the object's boundary and shape, resulting in a slightly coarse representation.

\begin{table}[!t]
\caption{Performance Metrics of the Deconvolution Experiment}\label{table:deconvolution}
\centering
\begin{tabular}{|c||c|c|c|c|c|}
\hline
 Method &Unknowns    &PSNR   &SNR    &SSIM  &MSE\\
 \hline
  \palentir{} &  \textbf{1125}  &  \textbf{32.7} & \textbf{24.4} & \textbf{53.4e-02}& \textbf{53.6e-05}\\
    \hline
 TV&   65536  &31.2 & 22.9 & {52.6e-02}& 75.8e-05\\
 \hline
\end{tabular}
\end{table}

\subsection{X-ray Computed Tomography}
\label{sec:experiment_7}

X-ray computed tomography (CT) is typically well-approximated as a linear problem and is one of the most common and well-known methods for medical imaging applications.
The forward model again takes the form of a matrix vector product, i.e., $\mathcal{M}(\bfpp) = \Rmb \bfpp$ where $\Rmb$ denotes a (discrete) Radon transform, and the vector $\mathcal{M}(\bfpp)$ denotes the vectorized form of the sinogram data. In this paper, we explore the performance of \palentir{} on both two and three dimensional CT problems, focusing specifically on limited-angle angle and/or sparse view reconstructions.  

\begin{figure}[!ht]
     \centering
     \subfloat[]
        {\includegraphics[width=0.32\textwidth]{tomo2/cheese_GT.pdf}
        \label{fig:cheese_GT}}
     \centering
     \subfloat[]
        {\includegraphics[width=0.32\textwidth]{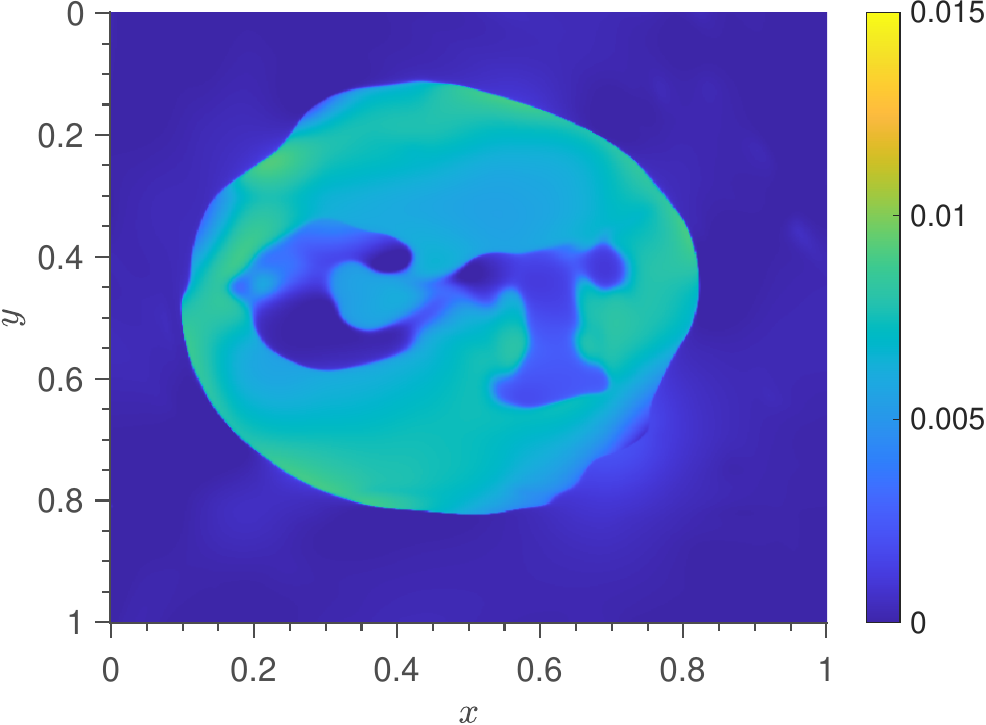}
        \label{fig:cheese_palentir}}
                 \centering
        \subfloat[]
        {\includegraphics[width=0.32\textwidth]{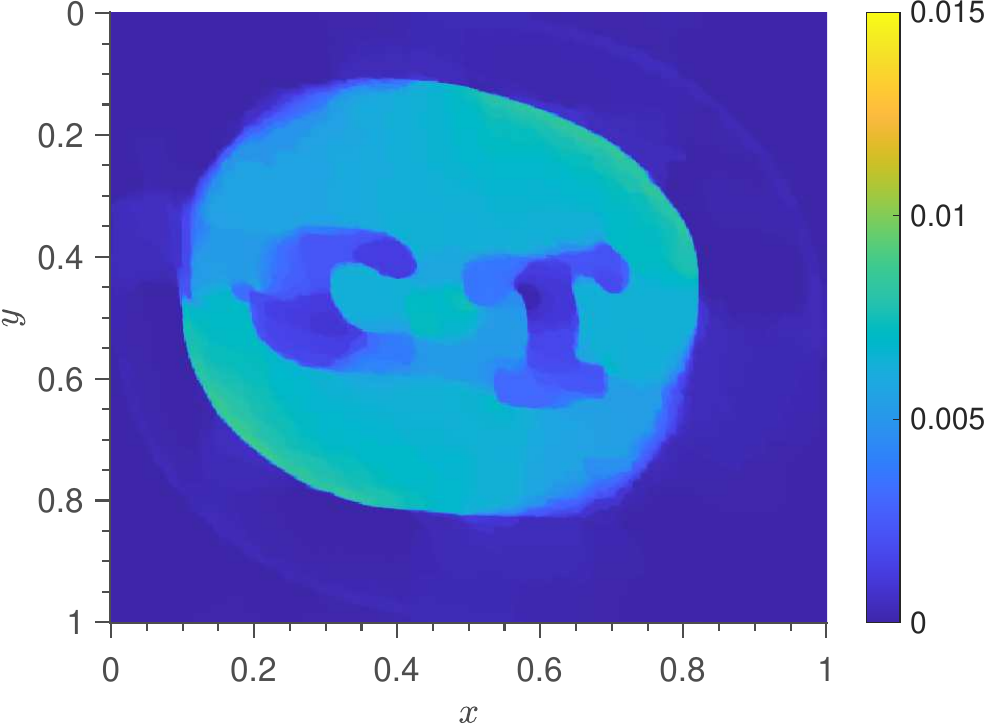}
        \label{fig:cheese_TV}}
     \centering \\
     \subfloat[]
        {\includegraphics[width=0.32\textwidth]{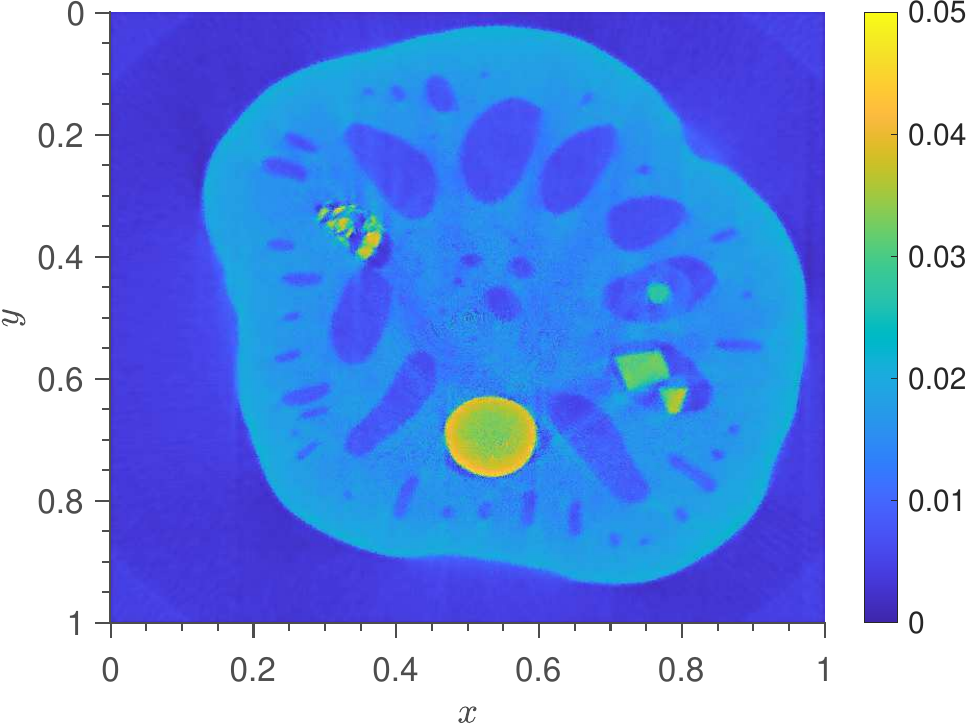}
        \label{fig:lotus_GT}}
                 \centering
                 \subfloat[]
        {\includegraphics[width=0.32\textwidth]{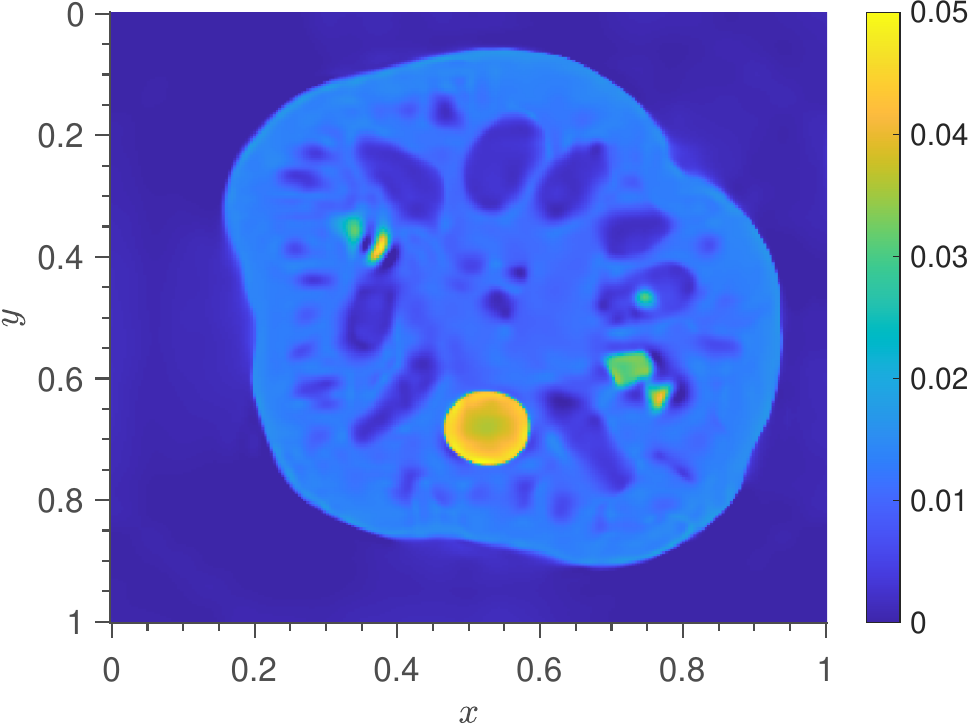}
        \label{fig:lotus_palentir}}
                 \centering
        \subfloat[]
        {\includegraphics[width=0.32\textwidth]{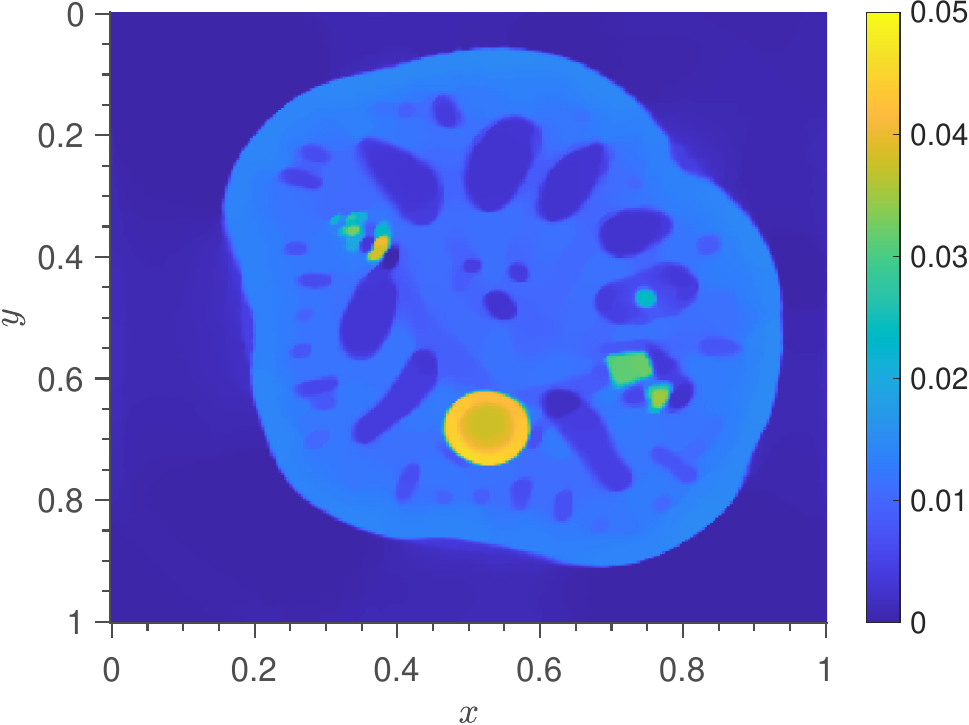}
        \label{fig:lotus_TV}}
                 \centering
\caption{2D Tomographic X-ray experiments: (a) Reference image of the carved cheese produced through high-resolution FBP reconstruction, computed from the full angle 360-projection sinogram. (b) \palentir{} limited angle and limited view reconstruction and (c) TV
limited angle and limited view
reconstruction of the carved cheese. (d) Reference image of the lotus root filled with attenuating objects, produced through high-resolution FBP reconstruction computed from the 360-projection sinogram. (e) Limited view \palentir{} reconstruction and (f) 
limited view TV reconstruction of the lotus root and the attenuating objects.}
\label{fig:tomo2}
\end{figure}

For the first 2D CT experiment, we use an open-access data set of tomographic X-ray data of a carved cheese \cite{bubba_cheese}. We note that, no approximation regarding the noise-level for this data set is given in \cite{bubba_cheese}. The reference image, generated by the given high-resolution filtered back-projection (FBP) reconstruction, is the $2000 \times 2000$ pixels image shown in Figure \ref{fig:cheese_GT}. In subsection \ref{sec:palentir_3}, we considered the sparse view problem with 15 projections spanning the full 360 degree circle. Now, we consider a limited angle, limited view problem with 15 projections spanning the range $1$\degree – $90$\degree. For this and the next CT experiments, there is a scaling and location difference between the provided high resolution filtered back-projection (FBP) and the CT data (which can be recognized from the recoveries); hence, we only provide a visual comparison of the TV and \palentir{} results for both 2D CT experiments. The $512 \times 512$ pixels \palentir{} and TV reconstructions of the carved cheese are shown in Figures \ref{fig:cheese_palentir} and \ref{fig:cheese_TV}. The limited angle, limited view nature of the data exacerbates the difficulty of the reconstruction. For the TV recovery, we see X-Ray artifacts, resulting in blurred boundaries of the cheese and carved letters.  Unintended artifacts in the background are also evident. The \palentir{} recovery of the 262144 pixels image of the carved cheese is produced using a grid of only $8 \times 8$ basis functions. Consequently, the number of unknowns is almost $92\%$ smaller than the number of unknowns required to solve for in a pixel-based method. The efficacy of the \palentir{} reconstruction is 
heavily influenced by the limitations of the 
available data, evident in imprecise boundaries at the top-left and bottom-right regions, 
similar to the errors in the TV recovery in these areas. However, in the remaining parts, \palentir{} demonstrates superior clarity and sharpness in both outer boundaries and intricate details around and within carved letters compared with the TV recovery.

For the second 2D experiment, we use an open-access data set of tomographic X-ray data of a lotus root filled with attenuating objects \cite{bubba16}. Again, an approximate noise-level is not given in \cite{bubba16}. The lotus root, akin in texture to a potato, features an array of holes of varying sizes; as it is primarily composed of starch, its structural characteristics make it an ideal candidate for stuffing with diverse objects. The holes contain the following four objects, each in a distinct hole within the lotus root: a pencil, a piece of chalk, three rectangular ceramic pieces, and several match-heads. Consequently, the lotus root, when filled with these objects, offers a diverse array of structures characterized by differing shapes, sizes, contrasts, and, notably, attenuations. This renders it a compelling subject for typical sparse-data CT applications \cite{bubba16}. The  $1500 \times 1500$ pixels reference image, shown in Figure \ref{fig:lotus_GT}, is generated from the high-resolution FBP reconstruction computed from the 360-projection sinogram provided in the dataset \cite{bubba16}. We consider a sparse view problem with 120 projections over the full range of angles. Figures \ref{fig:lotus_palentir} and \ref{fig:lotus_TV} show the $256 \times 256$ pixels \palentir{} and TV reconstructions. The TV method requires a considerably larger number of unknowns, specifically 65536, for the reconstruction of a $256 \times 256$ image than our \palentir{} approach, which utilizes a set of $25 \times 25$ ABFs, amounting to only 3125 parameters. Remarkably, both methods demonstrate 
a good recovery of fine image details. Considering the recovery of the lotus root itself, the TV method yields a smoother reconstruction compared with the \palentir{} recovery. Both reconstructions exhibit a smooth background, free from X-ray artifacts and effectively represent all four objects placed within the lotus root's hollows, even the smaller ones. Finally, it is worth highlighting that \palentir{} achieves the recovery of a fairly complex, multi-contrast, image from real X-ray data using
a single level-set function; a feature that distinguishes it from existing level-set methods. Yet, as a PaLS method, \palentir{} captures an enormous amount of detail that is even hard to recognize in the reference image.

\begin{figure}[!ht]
     \centering
     \subfloat[]
        {\includegraphics[width=0.49\textwidth]{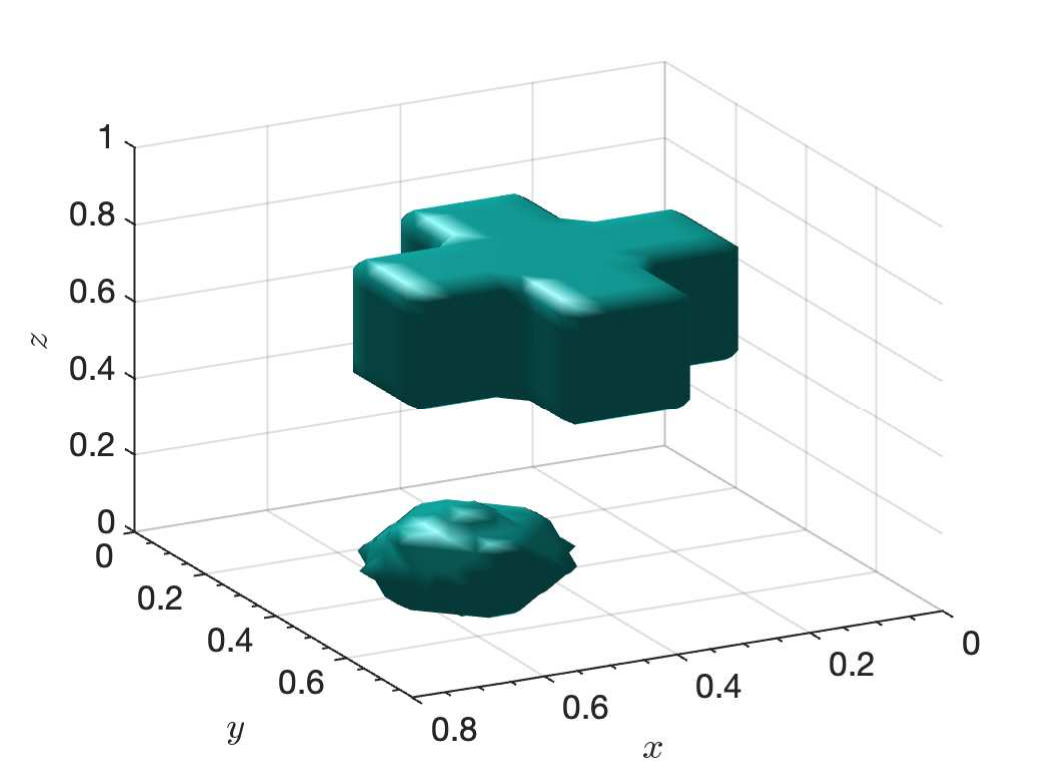}
        \label{fig:3d_true}}
     \centering
     \subfloat[]
        {\includegraphics[width=0.49\textwidth]{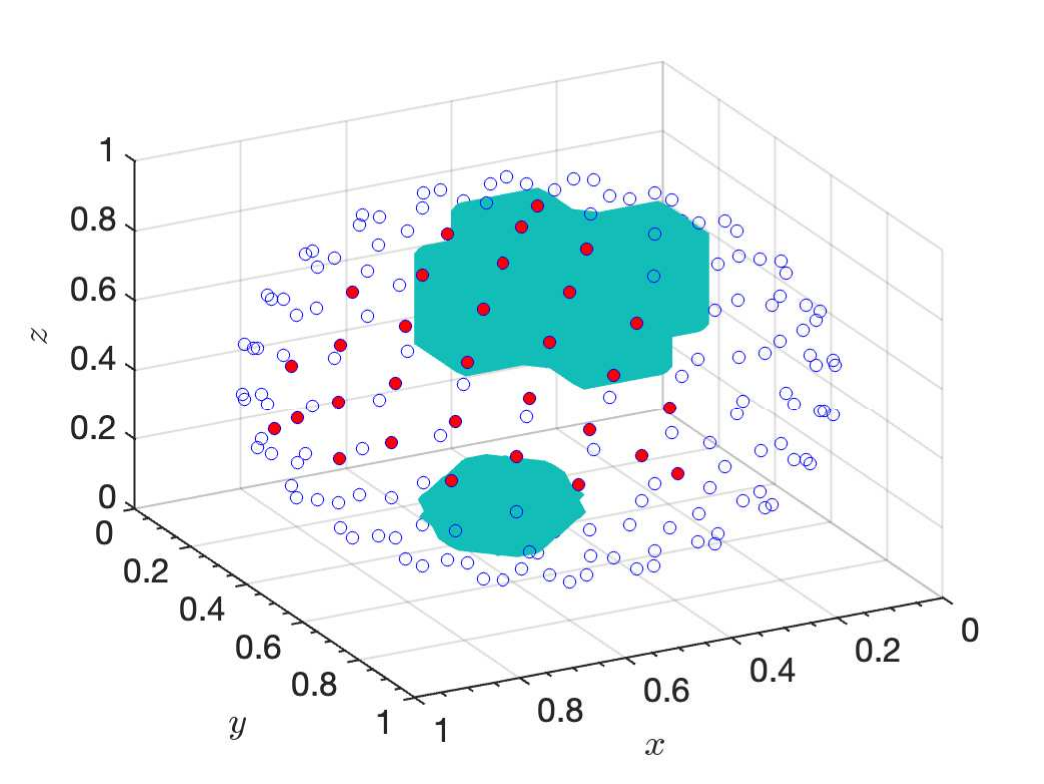}
        \label{fig:3d_projection}}
                 \centering \\
                 \subfloat[]
        {\includegraphics[width=0.49\textwidth]{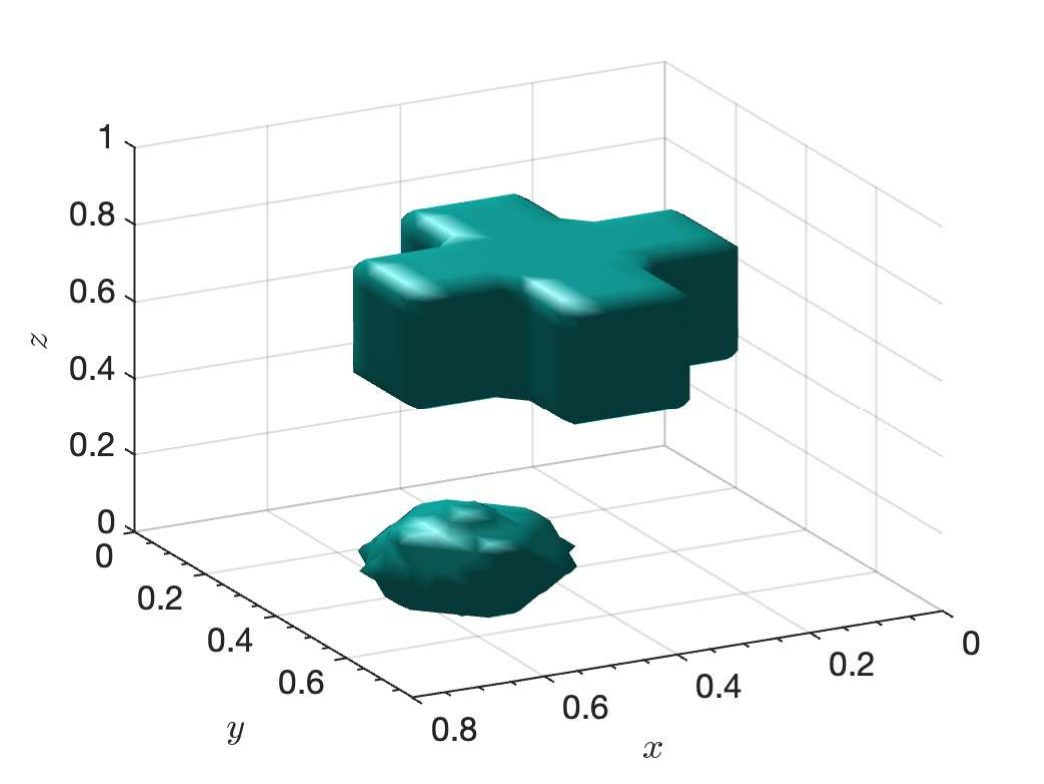}
        \label{fig:3d_palentir}}
                 \centering
     \subfloat[]
        {\includegraphics[width=0.49\textwidth]{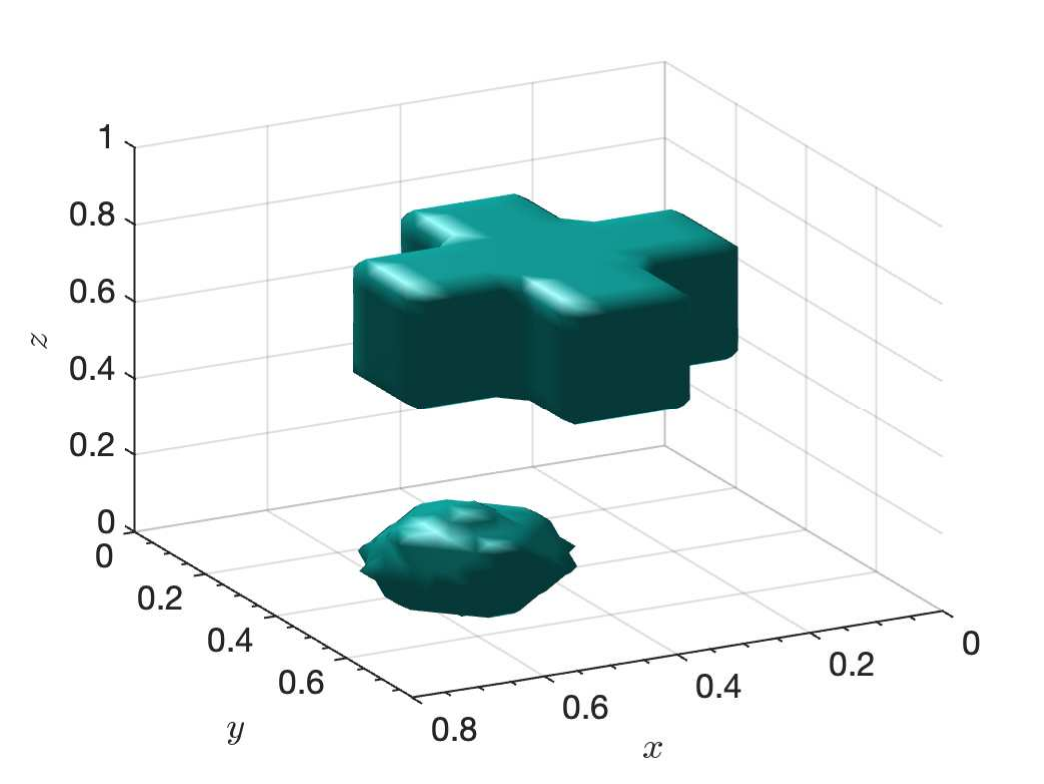}
        \label{fig:3d_TV}}
                 \centering
\caption{Results of 3D limited view parallel beam tomography experiment. (a) Original image of the 3D objects. (b) The distribution of projection plane centers. The 3D objects are represented in green, the blue points correspond to samples from a Lebedev sphere centered at the origin. The red points represent the selected samples from the sphere, chosen as the projection plane centers. (c) \palentir{} and (d) TV reconstructions. 
}
\label{fig:3d_figures}
\end{figure}

We now consider a 3D limited view parallel beam tomography experiment. We used \texttt{TVReg} \cite{hansen} to 
construct the experimental setup. The true 3D image is of size $27 \times 27 \times 27$, shown in Figure \ref{fig:3d_true}.
The input data is formed from 31 projections, each having a plane center located on the surface of the $1/8$ of a Lebedev sphere \cite{lebedev76}. Figure \ref{fig:3d_projection} illustrates the spatial distribution of projection plane centers for our 3D limited-angle view CT experiment. The objects of interest are represented in green, while the blue rings indicate the homogeneously spaced samples obtained from a Lebedev sphere centered at the origin. Notably, 31 of these blue rings are highlighted in red, indicating the selected projection plane centers for the experiment. These centers have been positioned to all reside within the same $1/8$ of the sphere, making the CT sparse and limited angle view.

The input data is corrupted with $1\%$ additive Gaussian white noise. In this experiment, \palentir{} is constructed using 343 basis functions centered on a $7 \times 7 \times 7$ grid resulting in 3087 unknowns. We conduct a performance comparison between \palentir{} and L2-Total Variation (TV) regularization. The reconstructions by \palentir{} and TV are presented in Figures \ref{fig:3d_palentir} and \ref{fig:3d_TV}, respectively, and the corresponding performance metrics are provided in Table \ref{table:3d}. In the context of a limited-angle view 3D CT problem, both the \palentir{} and the TV method exhibit remarkable performance. Both methods produce high values for PSNR, SNR, and small values for MSE, with the 
TV method performing marginally better.
The results for these metrics, all of which serve as critical metrics for assessing image quality and fidelity, 
emphasize the capacity for 
both methods to 
effectively reduce absolute error. Both methods attain the highest possible SSIM scores, a fact corroborated by the significant resemblance between the reconstructed images in Figures \ref{fig:3d_palentir} and \ref{fig:3d_TV} to the 
reference image shown in Figure \ref{fig:3d_true}. \palentir{} maintains its advantage by recovering a 19683 pixel 3D image using only 3087 parameters, 
without the need of explicit
regularization. 
\begin{table}[!ht]
\caption{Performance Metrics of the 3D CT Experiment}\label{table:3d}
\centering
\begin{tabular}{|c||c|c|c|c|c|}
 \hline
 Method &Unknowns    &PSNR   &SNR    &SSIM  &MSE\\
 \hline
 TV&   19683  & \textbf{61.0} & \textbf{48.0} & \textbf{1.0}& \textbf{80.3e-08}\\
  \hline
 \palentir{} &  \textbf{3087}  &  {60.1} & {47.2} & \textbf{1.0}& {96.9e-08}\\
 \hline
\end{tabular}
\end{table}


\subsection{Diffuse Optical Tomography} 
\label{sec:experiment_10}

Diffuse optical tomography (DOT) is a  non-invasive, low-cost alternative for breast and brain imaging compared with X-ray and MRI \cite{boas2001imaging}. In DOT, the tissue is illuminated with near-infrared light and the data, comprised of point measurements of diffused and partially absorbed photon fields, is collected external to the body. These measurements are used along with a mathematical model, typically a diffusion-absorption equation (posed in the frequency domain), to recover the optical absorption and (sometimes) scattering properties of the medium.  Here we use such a model of the form
\begin{equation}
-\nabla\cdot(D(\rmb)\nabla \eta(\rmb)) +
\mu(\rmb;{\bf p})\eta(\rmb) + \frac{\imath\omega}{\nu}\eta(\rmb)
 =  g(\rmb) ,
\end{equation} 
where $D(\rmb)$ represents the  (here) known scalar diffusion at a point $\rmb$, $\mu(\rmb;{\bf p})$, the quantity for which we invert, represents the absorption as a function of space and the parameter vector $\bf p$, $\omega$ represents the modulation frequency of the light source, and $\nu$ represents the speed of light in  the tissue. The sources are placed one grid point inside the medium \cite{arridge}, and the detectors are placed on the opposite boundary. On the boundaries where the sources and detectors are located, we have Robin boundary conditions; on the other boundaries, we use homogeneous Dirichlet boundary conditions, $\eta(\rmb) = 0$. For details, see \cite{StuEtAl2015,Aslan2019,antoulas2020interpolatory}. The recovered absorption coefficient (and sometimes also the diffusion/scattering) can be used to characterize the state of the tissue  \cite{boas2001imaging, fang2009combined}.

\begin{figure}[!ht]
     \centering
     \subfloat[]%
    {\includegraphics[width=0.32\textwidth]{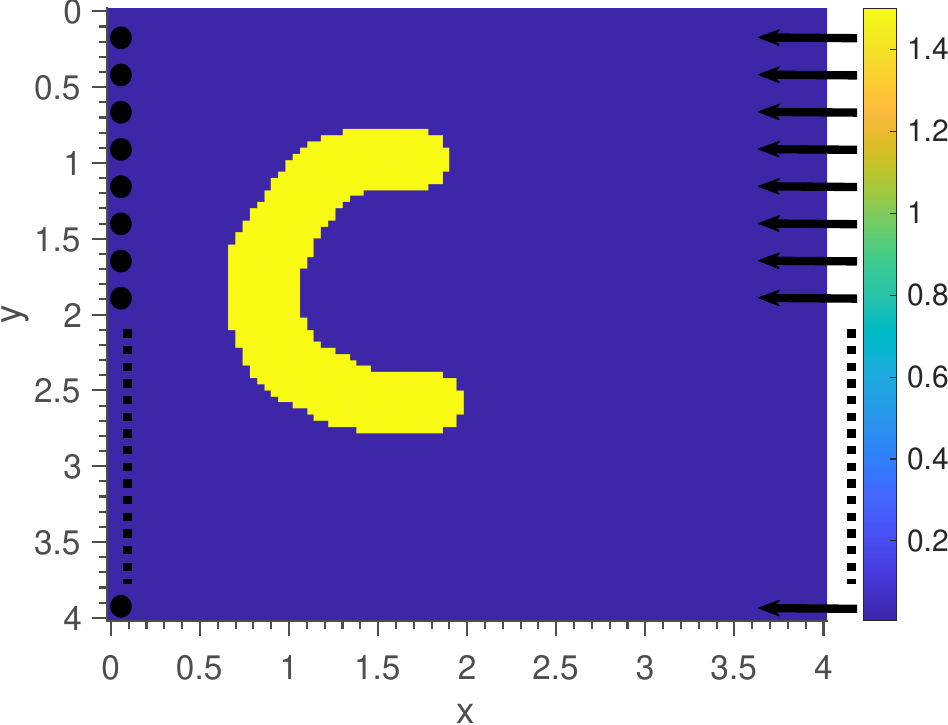}
    \label{fig:DOT_TrueAnomaly}}
        \centering
    \subfloat[]
        {\includegraphics[width=0.32\textwidth]{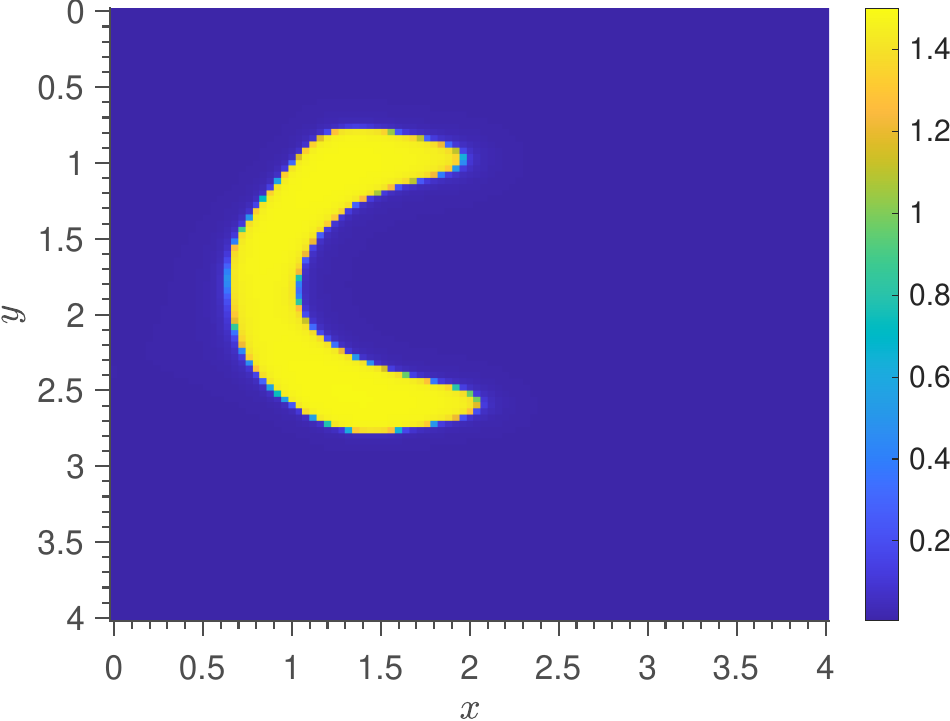}
        \label{fig:DOT_palentir_noiseless}}
                 \centering
    \subfloat[]
        {\includegraphics[width=0.32\textwidth]{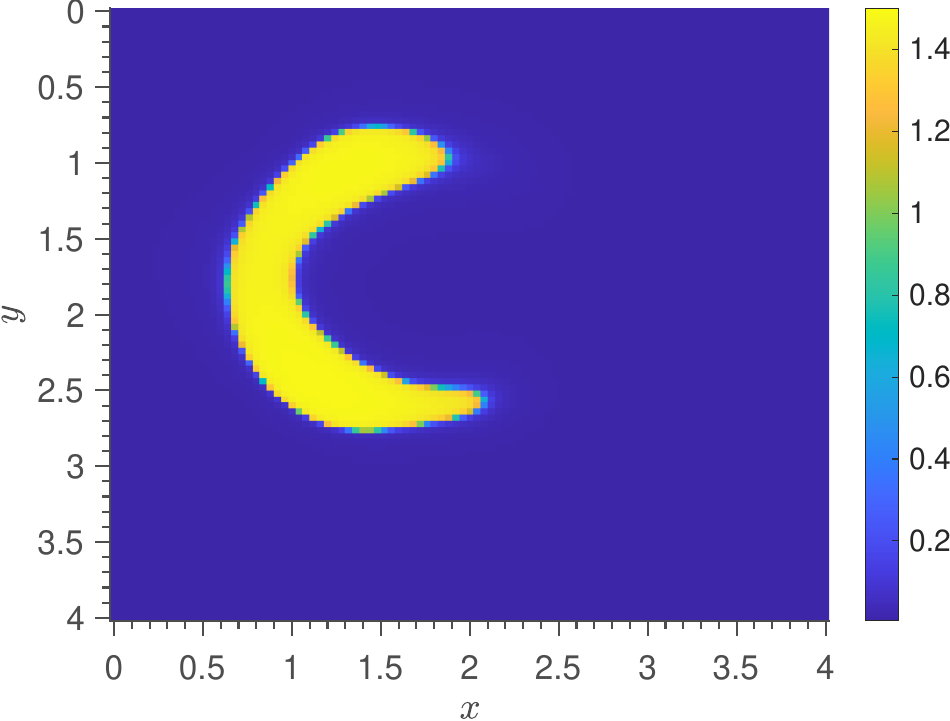}
        \label{fig:DOT_palentir_noise}}
                 \centering
        \caption{ The results of Diffuse Optical Tomography experiment are shown. (a) True anomaly with sources (right) and Detectors (left), (b) reconstruction of the absorption image without additive white noise (SSIM:0.92, MSE:1.46e-2, PSNR:21.87dB, SNR:10.75dB), (c) reconstruction of the absorption image with $1\%$ additive white noise (SSIM:0.89, MSE:1.92e-2, PSNR:20.69dB, SNR:9.57dB).}
\label{fig:DOT_Experiment}
\end{figure}

Following, e.g.,~\cite{larusson2012parametric}, we assume that the absorption coefficient can be modeled via \eqref{eq:newpals}.  As shown in Figure \ref{fig:DOT_TrueAnomaly}, we take the region to be imaged as a rectangle of size 4cm by 4cm, with $m_s = 32$ sources arrayed on the right side and $m_d = 32$ detectors on the left.  As before, we let ${\bf f}({\bf p})$ denote the discrete absorption image for a given parameter vector, that is, $[{\bf f}({\bf p})]_i = f({\bf r}_i,{\bf p})$ for any grid point ${\bf r}_i$. Assuming we collect data for all detectors when each source is active, we obtain a data vector $\bf d$ with $m = m_s \times m_d $ values for each modulation frequency $\omega$.

The input-output map from sources to detectors \cite{SAIBABA21} (also called the  \textit{transfer function}) as a function of $\bf p$ and $\omega$, is given by
\begin{equation} \label{eq:Psi}
  {\boldsymbol \Psi}(\bfpp;\omega) =  {\bf C}^T {\bf A}(\bfpp;\omega)^{-1} {\bf B} \quad \in \> \mathbb{R}^{m_d\times m_s},
\end{equation}
where ${\bf B} \in \mathbb{R}^{n \times m_s}$ represents $m_s$ sources, $n$ is the total number of voxels or grid points, 
${\bf A}({\bf f}({\bf p});\omega) 
\in \mathbb{R}^{n \times n}$ represents the discretization of the diffusion-absorption
equation, and 
${\bf C}^T \in \mathbb{R}^{ m_d \times n}$ simulates the measurement of outputs at $m_d$ detectors. So, ${\bf A}(\bfpp;\omega) {\bf X} = {\bf B} $ represents the discretized PDE that relates photon fluence/flux at grid points to the sources \cite{arridge}.  The DOT inverse problem is then specified by the forward mapping ${\cal M}(\fmb(\pemb);\boldsymbol{\omega}) = \text{vec}(\boldsymbol \Psi(\bfpp ;\boldsymbol{\omega} ))$, the vectorization of the transfer function outputs for a vector of frequencies 
$\boldsymbol{\omega} = [\omega_1, \ldots, \omega_{m_{f}}]$. Given a vector of
measured data (with additive noise) $\bf d$, we solve for $\bf p$ by minimizing
$\frac{1}{2} \| \text{vec}(\boldsymbol \Psi(\bfpp;
\boldsymbol{\omega} )) - \dmb \|_2^2$. Note that regularization is provided
implicitly by the
parameterization.

We present results in Figure
\ref{fig:DOT_Experiment}, with the true anomaly shown in Figure \ref{fig:DOT_TrueAnomaly}, and
we focus on the quality of the shape reconstruction using the new \palentir{} parameterization. The DOT problem is severely ill-posed, due to both the limited numbers of sources and detectors and the highly diffuse nature of the problem, so accurate reconstructions are not possible.  Indeed, this is why the shape-based parameterized reconstruction for problems with well-defined boundaries is important. To distinguish between reconstruction error due to noise as well ass ill-posedness and reconstruction error due to the nature of the problem, we provide a reconstruction with 1$\%$ additive Gaussian white noise (realistic) as well as one without noise. The latter serves to demonstrate the difficulty of the data limited problem which should be distinguished from the approximation quality of PaLEnTIR. For both \palentir{} recoveries, $11^2=121$ ABFs are used and contrast coefficients are chosen as minimum and maximum values of the ground-truth image. The performance metrics are given in the figure description for both recoveries. The reconstruction for the data without noise is shown in Figure \ref{fig:DOT_palentir_noiseless}. The model captures the structure quite well with SSIM score of 0.92, with slight imperfections due to the data limited nature of the problem rather than the PaLEnTIR model. When additive noise is added and we run the algorithm we observe the image in Figure \ref{fig:DOT_palentir_noise}. While the edges are slightly less clear in the noisy case, the reconstruction is remarkably good with SSIM score of 0.89, given the limitations of the data and the diffusive nature of the forward problem. It is important to note that there was no parameter to tune in generating this reconstruction.

\section{Conclusion and Future Work}
\label{sec:conclusion}
In this paper, we redefine the parametric level-set method to enhance the quality and functionality of the PaLS approach in the context of solving inverse imaging problems for images with (possibly multiple) piecewise constant contrasts.
As the most significant contribution of our work, the \palentir{} parameterization stands out as the only level-set approach, as far as we know, that employs only a single level-set function, irrespective of the number of contrasts or their values. Thus \palentir{} can represent multi-contrast scenes with very few parameters compared with traditional level-set based and pixel-based methods, and without any need to tune a regularization parameter. 
We investigate the qualitative performance of our new model 
and demonstrate that replacing RBFs in PaLS with ABFs expands the family of shapes that can be represented by a single basis function from circles to rotated ellipses. 
Relative to RBF PaLS, \palentir{} demonstrates a better reconstruction
capability by capturing significantly more details, by converging much faster, and 
by achieving significantly 
lower objective function values. We empirically demonstrate how \palentir{}  outperforms the RBF PaLS in terms of the condition number in both single and multiple basis functions cases. 
 
We demonstrate the utility of \palentir{} on numerical experiments over a range of  inverse problems, including 2D and 3D experiments, real data, and both non-linear and linear inverse problems.  In the examples where comparisons are made, \palentir{} reconstructions are about as good or somewhat better 
than those computed with competing reconstruction methods, whose parameters have been manually tuned to minimize the MSE, which is not
possible in practice.  
The built-in regularization through the parameterization allows, even in the non-linear, severely data-limited DOT problem, for high-quality reconstructions.  
Multiple contrasts in varying shapes and configurations are readily picked up by our new approach without any hand-tuning or regularization parameter selection. In summary, we have demonstrated the power of \palentir{} for reconstruction, denoising, and restoration for the class of 
piecewise constant images.       

The  performance exhibited by \palentir{} opens up prospects for future research in several key areas. Our forthcoming work will comprehensively explore the application of \palentir{} to three-dimensional domains, subjecting the model to diverse nonlinear and ill-posed problems. While the current use of homogeneously distributed basis functions has proven effective, we envision further enhancements by incorporating an adaptive refinement strategy, akin to strategies proposed in \cite{baussard2004adaptiveA, baussard2004adaptiveB}, which place additional basis functions in regions of higher geometric complexity. 
Additionally, our research agenda includes a focused investigation into the application of \palentir{} for uncertainty quantification. This entails the development of prior models for \palentir{} parameters based on object priors, with a specific emphasis on leveraging these priors to quantify accuracy in tasks such as object localization and characterization throughout the recovery process.

\appendix
\section{Derivatives of the 2D \palentir{} parameterization with respect to model parameters}
\label{sec:derivatives}

The \texttt{TREGS} algorithm necessitates a function for computing the Jacobian matrix. In this section, we present the derivations essential for calculating the Jacobian of the \palentir{} model, denoted as $f(\rmb;\pemb)$. Since $\pemb$ encompasses both anisotropic basis function (ABF) parameters and contrast parameters, we represent these parameters with two separate vectors: $\hat{\pemb}$ for the ABF parameters and $\mathbf{p}_c$ for the contrast parameters. The \palentir{} model, now expressed as $f(\rmb;\hat{\pemb}, \mathbf{p}_c)$, is defined as:
\begin{equation}
f(\rmb;\hat{\pemb};\mathbf{p}_c)  =
    C_{\Delta}(\rmb;\mathbf{p}_c) T\left(\phi(\rmb;\hat{\pemb}) \right)+ C_{L} (\rmb;\mathbf{p}_c).
\end{equation}
 Here $C_{\Delta}(\rmb;\mathbf{p}_c)=C_{H} (\mathbf{\rmb};\mathbf{p}_c)-C_{L} (\rmb;\mathbf{p}_c)$. We begin by examining the derivatives with respect to $\hat{\pemb}$. The derivative of $f(\rmb;\hat{\pemb};\mathbf{p})$ with respect to a parameter $ \hat{p}_i $, that is the $i^{th}$ parameter in $3N$-vector $\hat{\pemb}$, is expressed using the chain rule as:
\begin{equation}
    \label{eq:chain}
    \frac{\partial f(\rmb;\hat{\pemb};\mathbf{p}_c)}{\partial \hat{p}_i}=C_{\Delta}(\rmb;\mathbf{p}_c)T'(\phi(\rmb;\hat{\pemb}))
    \frac{\partial \phi(\rmb;\hat{\pemb})}{\partial p_i}.
\end{equation}
The first factor on the right side in \eqref{eq:chain} is the difference between the contrast coefficients. The second term is the derivative of the transition function. The transition function for a $c$-level-set, $T(x)$,  is defined as:
\begin{equation}
    \label{eq:app_transition}
    T\left( x \right)= \frac{1}{2} \left[ 1+ \frac{2}{\pi} \text{tan}^{-1}\left(\frac{\pi (x-c)}{w}\right) \right],
\end{equation}
and the derivative of \eqref{eq:app_transition} is given by:
\begin{equation}
    T'(\phi(\rmb;\hat{\pemb}))=\frac{1}{w}\left[\frac{1}{1+(\frac{\pi x}{w})^2}\right].
\end{equation}

The third factor in \eqref{eq:chain} is the derivative of the parametric level-set function, $\phi$, with respect to model parameters. $\phi$ is defined in \eqref{eq:newphi}. The expressions for the derivatives of \eqref{eq:newphi} with respect to each parameter are provided below, where elements $\left [ \alpha_j; \beta_j; \gamma_j  \right ]$ correspond to the parameters of the $j^{th}$ ABF for  $j=1,2 \dots N$:
\begin{align}
    \frac{\partial \phi(\rmb;\hat{\pemb})}{\partial \alpha_j}  &= \frac{\partial \sigma_h(\alpha_j)}{\partial \alpha_j}\psi \left( \Rmb_j (\rmb-\boldsymbol{\chi}_{j})\right), \\
    \frac{\partial \phi(\rmb;\hat{\pemb})}{\partial \beta_j} & = \sigma_h \left(\alpha_j \right) \frac{\partial \psi \left( \Rmb_j (\rmb-\boldsymbol{\chi}_{j})\right)}{\partial \beta_j},
    \\
    \frac{\partial \phi(\rmb;\hat{\pemb})}{\partial \gamma_j} & = \sigma_h \left(\alpha_j \right) \frac{\partial \psi \left( \Rmb_j (\rmb-\boldsymbol{\chi}_{j})\right)}{\partial \gamma_j}.
\end{align}
The derivatives with respect to parameters in the equations above are defined as:
\begin{align}
    \frac{\partial \sigma_h(\alpha_j)}{\partial \alpha_j} &= \frac{1}{2}\text{sech}^2(\frac{\alpha_j}{2}), \\
   \frac{\partial \psi \left( \Rmb_j \rmb_{\Delta}\right)}{\partial \beta_j} & =
   (-2\mu^2) \psi \left( \Rmb_j \rmb_{\Delta}\right)\rmb_{\Delta}^T 
    \begin{bmatrix}
                    e^{2\beta} & 0\\
                    \gamma e^\beta & -e^{-\beta} 
                \end{bmatrix}
   \rmb_{\Delta},
   \\
   \frac{\partial \psi \left( \Rmb_j \rmb_{\Delta}\right)}{\partial \gamma_j} & =
   (-\mu^2)\psi \left( \Rmb_j \rmb_{\Delta}\right) \rmb_{\Delta}^T 
    \begin{bmatrix}
                    0 & e^\beta\\
                    e^\beta & 2\gamma 
                \end{bmatrix}
   \rmb_{\Delta}.
\end{align}
where $\rmb_{\Delta}=(\rmb-\boldsymbol{\chi}_{j})$. Next, the derivative of $f(\rmb;\hat{\pemb};\mathbf{p}_c)$ with respect to $\mathbf{p}_c$ is defined as
\begin{equation}
    \label{eq:pcchain}
    \frac{\partial f(\rmb;\hat{\pemb};\mathbf{p}_c)}{\partial pc_k}=C_{\Delta}(\rmb;\bm{e_k}) T\left(\phi(\rmb;\hat{\pemb}) \right)+ C_{L} (\rmb;\bm{e_k}),
\end{equation}
where $\bm{e_k}\in \R^{2N}$ is a $2N-$vector with it's $k^{th}$ element equal to $1$ and all other elements equal to $0$, for  $k=1,2 \dots 2N$.
\section*{References}
\bibliography{palentir} 

\end{document}